\documentclass[12pt]{amsart}
\usepackage{latexsym,amsmath, amscd, amssymb, amsthm, euscript, mathrsfs}
\usepackage[all]{xy}
\newtheorem{thm}[subsection]{Theorem}
\newtheorem{thm/def}[subsection]{Theorem/Definition}
\newtheorem{cor}[subsection]{Corollary}
\newtheorem{lem}[subsection]{Lemma}

\newtheorem{prop}[subsection]{Proposition}
\newtheorem{problem}[subsection]{Problem}
\theoremstyle{definition}
\newtheorem{notation}[subsection]{Notation}

\theoremstyle{definition}

\theoremstyle{definition}

\newtheorem{rem}[subsection]{Remark}

\numberwithin{equation}{subsection}

\newtheorem{pg}[subsection]{}

\newcommand{\Q}{\mathbb{Q}}
\renewcommand{\P}{\mathbb{P}}
\newcommand{\A}{\mathbb{A}}

\newcommand{\Z}{\mathbb{Z}}

\newcommand{\Sp}{\text{\rm Spec}}

\newcommand{\mls}{\ms} 

\DeclareMathOperator{\spec}{Spec}
\newcommand{\ms}{\mathscr}
\newcommand{\simto}{\stackrel{\sim}{\to}}

\newcommand{\surj}{\twoheadrightarrow}
\newcommand{\mf}{\mathfrak}
\newcommand{\widebar}{\overline}

\newcommand{\m}{\boldsymbol{\mu}}

\title{Generators and relations for the \'etale fundamental group}

\author{Max Lieblich}
\address{Fine Hall, Washington Road, Princeton NJ 08544-1000}

\author{Martin Olsson}
\address{Dept. of Mathematics, University of California, Berkeley, CA 94720-3840}

\begin{document}

\begin{abstract}If $C$ is a smooth curve over an algebraically closed
  field $k$ of characteristic exponent $p$, then the structure of the
  maximal prime-to-$p$ quotient of the \'etale fundamental group is
  known by analytic methods.  In this paper, we discuss the properties
  of the fundamental group that can be deduced by purely algebraic
  techniques.  We describe a general reduction from an arbitrary curve
  to the projective line minus three points, and show what can be
  proven unconditionally about the maximal pro-nilpotent and
  pro-solvable quotients of the prime-to-$p$ fundamental group.
  Included is an appendix which treats the tame fundamental group from
  a stack-theoretic perspective.
\end{abstract}

\maketitle

\section{Introduction}

Let $k$ be an algebraically closed field of characteristic exponent
$p$.  For a connected scheme $X/k$ with a geometric
point $\bar x\rightarrow X$, let $\pi _1'(X, \bar x)$ denote the
maximal prime-to-$p$ quotient of the \'etale fundamental group of $X$.

The motivation for this note is the following problem, posed in
\cite[X.2.7]{SGA1}.

\begin{problem}\label{hardproblem}
  Let $(C, \bar x\rightarrow C)/k$ be a connected smooth proper
  pointed curve of genus $g$ over $k$.  Give an algebraic proof of the
  fact (which follows from the comparison between the \'etale
  fundamental group and the topological fundamental group over
  $\mathbb{C}$) that the group $\pi _1'(C, \bar x)$ is the maximal
  prime-to-$p$ quotient of the profinite completion of the free group
  on $2g$ generators $\gamma _1, \sigma _1, \gamma _2, \sigma _2,
  \dots , \gamma _g, \sigma _g$ modulo the relation
  \begin{equation}\label{relation} [\gamma _1, \sigma _1][\gamma _2,
    \sigma _2]\cdots [\gamma _g, \sigma _g] = 1.
  \end{equation}
\end{problem}

This appears to be a deep problem.  For example, a solution would
likely shed light on the inverse Galois problem.  

There is also a version of Problem \ref{hardproblem} for punctured curves.

\begin{problem}\label{hardproblembis}
  Let $C/k$ be a smooth complete connected curve, $p_1, \dots , p_n\in C(k)$
  distinct points ($n\geq 1$), and let $\widetilde{C}
  =C-\{p_i\}_{i=1}^n$.  Find an algebraic proof that for any geometric
  point $\bar x\rightarrow \widetilde{C}$ the group $\pi
  _1'(\widetilde{C}, \bar x)$ is isomorphic to the maximal
  prime-to-$p$ quotient of the profinite completion of the free group
  on $2g+n-1$ generators.
\end{problem}

It is apparently known to experts that in order to solve Problem
\ref{hardproblem} it suffices to solve Problem \ref{hardproblembis}
for $\mathbb{P}^1-\{0,1,\infty \}$ (in particular, Fujiwara has
apparently obtained many of the results of this paper using the log
\'etale fundamental group \cite{F}).  Our aim in this paper is to
write out some aspects of this bit of folklore, and to explain some
properties of the fundamental group that can be proven using algebraic
methods.

In order not to confuse the reader, from now on all stated results
will be proven using only algebraic methods (so the content of the
following theorems is that they can be proven purely algebraically).
Write $F_{n}$ for the free profinite group on $n$ generators and write
$S_g$ for the quotient of $F_{2g}$ with generators
$\gamma_1,\sigma_1,\dots,\gamma_g,\sigma_g$ by the closed normal
subgroup generated by $[\gamma_1,\sigma_1]\cdots[\gamma_g,\sigma_g]$.
In what follows, the superscript $t$ indicates the tame fundamental
group, while the superscript $'$ indicates the prime-to-$p$ quotient
of the fundamental group.  Similarly, given a profinite group $G$ and
a class of finite groups $\ms C$ satisfying the conditions in
Paragraph \ref{2.5}, we write $G^{\ms C}$ to indicate the maximal
pro-$\ms C$ quotient.  (The reader is referred to
\S\ref{sec:mc-c-groups} for definitions and notation.)  For example,
the superscript $\text{sol}$ (resp.\ $\text{nil}$) indicates the
maximal pro-solvable (resp.\ maximal pro-nilpotent) quotient.

\begin{thm}\label{theorem} Let the notation be as in
  Problem \ref{hardproblem}.  
  \begin{enumerate}
  \item[$(i)$] There is an isomorphism $\pi_1'(C,\bar
    x)^{\text{\rm nil}}\simeq (S_g')^{\text{\rm nil}}$.
  \item[$(ii)$] Suppose the groups in $\mls C$ have orders prime to
    $p$.  If 
    $\pi_1^{\ms C}(\P^1\setminus\{0,1,\infty\},\bar x)$ is
    topologically freely generated by a pair of ramification elements
    at $0$ and $1$ (see \S\ref{section8}), then the maximal pro-$\ms
    C$ quotient $\pi _1^{\ms C}(C, \bar x)$ is
    isomorphic to the maximal pro-$\ms C$ quotient $F_{2g}^{\ms C}$
    modulo the closed normal subgroup generated by a 
    non-trivial element lying in the closure of the commutator
    subgroup $[F_{2g}^{\ms C},F_{2g}^{\ms C}]$.
  \end{enumerate}
\end{thm}

Regarding Problem \ref{hardproblembis} we have the following.

\begin{thm}\label{theorembis} Let the notation be as in
  Problem \ref{hardproblembis}, and let $\ms C$ be a class of finite
  groups of prime-to-$p$ orders satisfying the conditions in Paragraph
  \ref{2.5}.
  \begin{enumerate}
  \item[$(i)$] If $\pi_1^{\ms C}(\P^1\setminus\{0,1,\infty\},\bar x)$
    is topologically freely generated by a pair of ramification
    elements at $0$ and $1$, then 
    $\pi _1^{\ms C}(\widetilde{C}, \bar x)$ is isomorphic to
    $(F_{2g+n-1})^{\ms C}$.
  \item[$(ii)$] In general 
    $\pi_1'(\widetilde C,\bar x)^{\text{\rm sol}}$ is topologically
    finitely generated.
  
 \item [$(iii)$] There is an isomorphism $\pi_1'(\widetilde C,\bar 
x)^{\text{\rm nil}}\simeq (F_{2g+n-1}')^{\text{\rm nil}}.$
  \end{enumerate}
\end{thm}

We remark that Theorem \ref{theorembis}$(\operatorname{ii})$ is clearly false without the
prime-to-characteristic hypothesis: it is easy to see by elementary
calculations in \'etale cohomology that the abelianization is not a free
$\Z$-module, as its $p$-part is too big by Abhyankar's conjecture
(proven by Raynaud \cite{raynaud} for the punctured projective line
and Harbater \cite{harbater} in general).

Using these results we can also prove using only algebraic techniques
the following finiteness results.

\begin{thm}\label{higherdim} Let $X/k$ be an irreducible scheme of
  finite type and $\bar x\rightarrow X$ a geometric point.  Then  $\pi _1^{\prime }(X, \bar x)^{\text{\rm sol}}$ 
 is topologically finitely
  generated.
\end{thm}

\begin{rem}
  A priori, it may seem that $\pi_1^t(X,\bar x)^{\text{\rm sol}}$
  depends upon a chosen compactification of $X$.  We refer the reader
  to Paragraphs \ref{firstpara} and \ref{secondpara} for basic facts
  about the tame fundamental group, and to Lemma
  \ref{L:A2.1} for the fact that the notion is
  independent of the choice of compactification (and in fact does not
  require a regular compactification to exist).
\end{rem}

This also gives a purely algebraic proof of an open version of the
Lang-Serre theorem \cite[X.2.12]{SGA1} for solvable groups.

\begin{cor}\label{C:lang-serre} Let $X/k$ be a scheme of finite type,
  and let $G$ be a finite solvable group of order prime to $p$.  Then the number of tamely
  ramified Galois covers of $X$ with group $G$ is finite.
\end{cor}

In X.1.10 of \cite{SGA1}, the reader will find wildly-ramified
counterexamples to Corollary \ref{C:lang-serre} for $X=\A^1_k$ and $G=\Z/p\Z$
(in the form of a continuous family of Artin-Schreier coverings).

The key to all of the results in this paper is a general d\'evissage
which enables one to deduce results for arbitrary curves from results
for $\mathbb{P}^1-\{0, 1, \infty \}$.  This is described in
\S\ref{section8}.

\begin{rem} With notation as in Problem \ref{hardproblembis}, the main
  result of Borne and Emsalem in \cite{Borne} gives an abstract
  isomorphism between the prime-to-$p$ prosolvable quotient
  $F_{2g+n-1}^{'\text{solv}}$ and $\pi _1^{\prime }(\widetilde C, \bar
  x)^{\text{solv}}$.  Their method, however, does not show that the
  generators can be chosen to be ramification elements in the sense of
  \S \ref{section8}.
\end{rem}

\begin{notation}\label{notrem}
If $G$ is a profinite group and if $X\subset G$ is a subset, then the \emph{closed normal subgroup} of $G$ generated by $X$ is by definition the intersection of all closed normal subgroups of $G$ containing $X$.

If $g$ and $n$ are nonnegative integers, we denote by $\mls M_{g, n}$
the Artin stack of $n$-pointed proper smooth curves of genus $g$.

We write $\widehat {\Z}(1)$ for $\varprojlim _{N}\m _N(\overline \Q)$.
\end{notation}

\begin{pg}{\bf Acknowledgments.} The authors are grateful to Niels
  Borne, David Harbater, Nick Katz, Fabrice Orgogozo, Brian
  Osserman, Bjorn Poonen, Akio Tamagawa, and Ravi Vakil for helpful conversations, and to Brian Conrad, Tam\'as Szamuely, and the referee for very careful readings of an earlier version of this paper.
  Lieblich was partially supported by an NSF Postdoctoral Fellowship
  and Olsson was partially supported by NSF grant DMS-0714086 and an
  Alfred P. Sloan Fellowship. \end{pg}

\section{$\ms C$-groups}
\label{sec:mc-c-groups}
\begin{pg}\label{2.1}
  It is useful to have a formalism for dealing with various quotients
  of $\pi_1$ in a uniform manner.  We adopt a well-known method and
  include this section primarily for the sake of notation.

Let $\ms C$ be a full subcategory of the category of finite groups
such that 

\begin{enumerate}
\item If $G,H\in\ms C$ then $G\times H\in\ms C$.
\item If $H\subset G$ and $G\in\ms C$ then $H\in\ms C$.
\end{enumerate}

Given a profinite group $\mf G$, there is a \emph{maximal pro-$\ms C$
  quotient\/} $\mf G\to\mf G^{\ms C}$ which gives a left adjoint to
the natural inclusion of the category of pro-$\ms C$ groups
in the category of profinite groups.  In order to adhere to common
conventions, given a pointed scheme $(X,\ast)$ we will write $\pi_1^{\ms
  C}(X,\ast)$ for $\pi_1(X,\ast)^{\ms C}$.  The categories $\ms C$ which we
will encounter most are the following.

\begin{enumerate}
\item $\ms C={}'$, the category of finite groups of order prime to the
  characteristic exponent of the base field.
\item $\ms C=\text{\rm sol}$, the category of solvable finite groups.
\item $\ms C=\text{\rm nil}$, the category of nilpotent finite groups.
\item $\ms C=\ell$, the category of finite $\ell$-groups for a prime number $\ell $.
\end{enumerate}
Since every nilpotent group is solvable, there is a natural surjection
$\mf G^{\text{\rm sol}}\surj\mf G^{\text{\rm nil}}$ for any profinite
group $\mf G$.  
\end{pg}
\begin{pg}\label{P:exact}
Note also that the functor $\mf G\mapsto \mf G^{\mls C}$ is right
exact in the sense that if
$$
\mf G'\rightarrow \mf G\rightarrow \mf G^{\prime \prime }\rightarrow 1
$$
is an exact sequence of profinite groups, then the induced sequence
$$
  \mf G^{\prime \mls C}\rightarrow \mf G^{\mls C}\rightarrow \mf G^{\prime \prime \mls C}\rightarrow 1
$$
is also exact.  This follows from the universal property of the
morphism $\mf G\rightarrow \mf G^{\mls C}$.
\end{pg}
\begin{pg}\label{freedef}
As usual, given a set $X$, one can form the \emph{free pro-$\ms C$
  group\/} on $X$, denoted $F_X^{\ms C}$.  A concrete realization is
given by taking the pro-$\ms C$ completion of the ordinary free group
on $X$.  The cardinality of $X$ is the \emph{rank\/} of the free
group; when $X\cong\{1,\ldots,n\}$, we will write $F_n^{\ms C}$ for
$F_X^{\ms C}$.
\end{pg}

\begin{pg}\label{2.5} {\em We will assume for the remainder of this paper that $\ms C$
  satisfies the following conditions\/}.
\begin{enumerate}
\item If $G, H\in \mls C$ then $G\times H\in \mls C$.
\item If $G$ is in $\mls C$ then any subgroup of $G$ and any quotient
  group of $G$ 
  is in $\mls C$.
\item $\ms C$ is closed under extensions and contains a non-trivial group.
\end{enumerate}
\end{pg}

  The reader will observe that these conditions apply to all
of the classes of groups listed above, with the exception of the
nilpotent groups.  (In general, one must sacrifice certain nice
functorial properties in return for the strong structure theory of
nilpotent groups, which we touch upon below.)

Observe also that if $\mls C$ satisfies the above conditions, then
there exists a prime number $\ell $ such that $\mls C$ contains all
nilpotent $\ell $-groups.  Indeed (3) implies that there exists a
prime number $\ell $ such that $\mls C$ contains $\Z/\ell \Z$, and
then an inductive argument using (2) implies that any nilpotent $\ell
$-group is also in $\mls C$.

\begin{lem}\label{2.3b}
  Suppose 
$$1\to\mf G'\to\mf G\to H\to 1$$
is an exact sequence of profinite groups with $H$ in $\ms C$ (so in
particular $H$ is finite).  If $\mf G'\to{\mf G'}^{\ms C}$ is an
isomorphism then $\mf G\to{\mf G}^{\ms C}$ is an isomorphism.
\end{lem}
\begin{proof}
  It suffices to show that for any open normal subgroup $W\subset \mf
  G$ there exists an open normal subgroup $V$ of $\mf G$ contained in
  $W$ such that $\mf G/V$ is in $\ms C$.  Since $\mf G\to H$ is
  continuous, $V:=W\cap\mf G'$ is open and normal in $\mf G$. Since
  $\mf G'\to{\mf G'}^{\ms C}$ is an isomorphism, there exists an open
  normal subgroup $U$ of $\mf G'$ contained in $V$ such that the quotient
  $\mf G'/U$ is in $\mls C$.  By property (2) above, $\mf G'/V$, being
  a quotient of $\mf G'/U$, is in $\ms C$, whereupon, by property (3),
  $\mf G/V$ is in $\ms C$.
\end{proof}

\begin{lem}\label{L:exact} Let
\begin{equation}\label{exact1}
1\rightarrow \mf G'\rightarrow \mf G\rightarrow H\rightarrow 1
\end{equation}
be a short exact sequence of profinite groups, with $H$ an element of
$\ms C$.  Then the induced sequence of pro-$\mls C$-completions
\begin{equation}\label{exact2}
  1\rightarrow \mf G^{\prime \mls C}\rightarrow \mf G^{\mls C}\rightarrow H\rightarrow 1
\end{equation}
is exact.
\end{lem}
\begin{proof}
As mentioned in Paragraph \ref{P:exact}, the exactness of the sequence
$$
\mf G^{\prime \mls C}\rightarrow \mf G^{\mls C}\rightarrow H\rightarrow 1
$$
follows from the universal property of the pro-$\mls C$-completion functor.  Let $K'$ denote the kernel of the homomorphism $\mf G'\rightarrow \mf G^{\prime \mls C}$.

The universal property of the map $\mf G'\rightarrow \mf G^{\prime \mls C}$ implies that for any element $g\in \mf G$ there exists a unique homomorphism $\rho _g:\mf G^{\prime \mls C}\rightarrow \mf G^{\prime \mls C}$ such that the diagram
$$
\xymatrix{\mf G'\ar[d]\ar[r]^{h\mapsto ghg^{-1}}& \mf G'\ar[d]\\
\mf G^{\prime \mls C}\ar[r]^{\rho _g}& \mf G^{\prime \mls C}}
$$
commutes.  This implies that the kernel $K'$ is also normal as a
subgroup of $\mf G$.  Let $\widetilde {\mf G}$ denote the quotient
$\mf G/K'$ so that there is an exact sequence
$$
1\rightarrow \mf G^{\prime \mls C}\rightarrow \widetilde {\mf G}\rightarrow H\rightarrow 1.
$$
By Lemma \ref{2.3b}, the map $\widetilde {\mf G}\rightarrow \widetilde
{\mf G}^{\mls C}$ is an isomorphism, and therefore we obtain a
commutative diagram
$$
\xymatrix{
  &\mf G^{\prime \mls C}\ar[d]^{\text{id}}\ar[r]&\mf G^{\mls C}\ar[r]\ar[d]&H\ar[d]^{\text{id}}\ar[r]&1\\
  1\ar[r]& \mf G^{\prime \mls C}\ar[r]& \widetilde {\mf G}\ar[r]&
  H\ar[r]& 1.}
$$
It follows that $\mf G^{\prime \mls C}\rightarrow \mf G^{\mls C}$ is
injective (and in fact that $\mf G^{\mls C}\rightarrow \widetilde {\mf
  G}$ is an isomorphism).
\end{proof}

A fundamental result related to free pro-$\ms C$ groups (see Paragraph
\ref{freedef}) we will use is the pro-$\ms C$ version of the
Nielsen-Schreier Theorem (in a weak form).  We give the statement
here for the purposes of convenience and refer the reader to Theorem
5.4.4 of \cite{Wilson} for the full statement of the theorem and its
proof (which is purely algebraic; the reader following the reference
can use the results of \S{6.1} of \cite{Robinson} to provide the
necessary algebraic proof of the classical version of the theorem).

\begin{thm}[Pro-$\ms C$ Nielsen-Schreier, weak form]\label{T:N-S}
  If $F$ is a free pro-$\ms C$ group of finite rank $r_F$ and
  $H\subset F$ is an open subgroup of index $i$, then $H$ is a free
  pro-$\ms C$ group of rank $r_H=i(r_F-1)+1$.
\end{thm}

Another useful fact is the following ``Hopfian property'' of  profinite groups:

\begin{thm}[{\cite[2.5.2]{Ribes}}]\label{Hopf}
  Let $F$ be a topologically finitely generated profinite group and
  let $\varphi :F\rightarrow F$ be a surjective endomorphism.  Then
  $\varphi $ is an isomorphism.
\end{thm}

\begin{pg}
  Finally, for future purposes, we remind the reader of two basic
  facts about pro-nilpotent groups.
\end{pg}
\begin{lem}\label{L:pro-nilpotent-dumb}
  Given a profinite group $\mf G$, there is a natural topological
  isomorphism $\mf G^{\text{\rm nil}}\simto\prod_{\ell}\mf G^{\ell}$,
  where the right side is given the product topology.
\end{lem}
\begin{proof}
  This follows from (1) the fact that any $\ell$-group is nilpotent,
  and (2) the fact that any nilpotent group is naturally isomorphic to
  the product of its Sylow subgroups.  See Prop.\ 2.4.3 of
  \cite{Wilson} for more details.
\end{proof}

\begin{lem}
  Let $\mf N$ be a pro-nilpotent group.  Given generators
  $x_1^{(\ell)},\dots,x_n^{(\ell)}\in\mf N^{\ell}$ for each prime
  $\ell$, there are generators $x_1,\dots,x_n$ for $\mf N$ such that
  the image of $x_i$ in $\mf N^{\ell}$ is $x_i^{(\ell)}$ for each
  prime $\ell$.
\end{lem}
\begin{proof}
  By Lemma \ref{L:pro-nilpotent-dumb}, there is a natural isomorphism
  $\mf N\to\prod_{\ell}\mf N^{\ell}$.  Thus, the elements
  $x_i^{(\ell)}$ give rise to elements $x_1,\dots,x_n$ of $\mf N$.  To
  see that these generate $\mf N$, it suffices to show that they
  generate any finite discrete nilpotent quotient of $\mf N$.  Thus,
  it suffices to establish the following: if $\ell _1, \dots, \ell _m$
  is a set of pairwise distinct primes and if for each $i=1, \dots, m$
  we have an $\ell _i$-group $G_{\ell _i}$ with a given choice of
  generators $x_1^{(\ell _i)}, \dots, x_n^{(\ell _i)}$, then the
  elements $x_j:= \prod _{i=1}^mx_j^{(\ell _i)}$ generate the group
  $G:= G_{\ell _1}\times \cdots \times G_{\ell _m}$.  This is a
  trivial application of the Chinese Remainder Theorem which we leave
  to the reader.
\end{proof}

\section{Ramification maps}\label{section8}

\begin{pg} Let $C/k$ be a connected smooth (not necessarily proper)
  curve, let $p\in C(k)$ be a point, and let $C^\circ $ denote
  $C-\{p\}$.  Fix a geometric generic point
  \begin{equation}\label{4.1.1b}
\bar \eta :\Sp (\Omega )\rightarrow C.
\end{equation}
Let $\mls O_{C, p}^{\text{sh}}$ denote the strict henselization of the
local ring at $p$, and let $\Sp (\mls O_{C, p}^{\text{sh}})^\circ $
denote the complement of the closed point in $\Sp (\mls O_{C,
  p}^{\text{sh}})$.

The tame  fundamental group of $\Sp (\mls O_{C, p}^{\text{sh}})^\circ $ with
respect to any base point is canonically isomorphic to $\widehat
{\Z}'(1):= \varprojlim _{(N, \text{char}(k)) =1}\m _N(\widebar{\Q})$.  The map
\begin{equation}\label{4.1.4c}
\Sp (\mls O_{C, p}^{\text{sh}})^\circ \rightarrow C^\circ 
\end{equation}
therefore induces for every choice of dotted arrow $h$  filling in the diagram
\begin{equation}\label{4.1.5b}
\xymatrix{
& \Sp (\Omega )\ar@{-->}[ld]_-{h}\ar[rd]^-{\bar \eta }& \\
\Sp (\mls O_{C, p}^{\text{sh}})^\circ \ar[rr]&& C^\circ }
\end{equation}
 a homomorphism
\begin{equation}\label{4.1.3}
\xymatrix{
\widehat \Z(1)\ar@{->>}[r]& \widehat \Z'(1)\ar[r]& \pi _1^t(C^\circ , \bar \eta ).}
\end{equation}
This homomorphism is well-defined up to conjugation in the following sense: if $\gamma \in \pi _1^t(C^\circ , \bar \eta )$, then the composition of (\ref{4.1.3}) with conjugation by $\gamma $ is also induced by a choice of $h$ filling in (\ref{4.1.5b}). We denote this conjugacy class of maps
by $c_p$. In order to specify a specific element of this conjugacy
class of homomorphisms, it suffices to choose an embedding $\mls O_{C,
  p}^{\text{sh}}\hookrightarrow \Omega $ compatible with the given
embedding $\mls O_{C, p}\hookrightarrow \Omega $.  We refer to the
homomorphism
$$
\rho :\widehat {\Z}(1)\rightarrow \pi _1^t(C^\circ , \bar \eta )
$$
obtained by making such a choice as the \emph{ramification map at $p$}.
\end{pg}

\begin{pg} Let $K\subset \pi _1^t(C^\circ , \bar \eta )$ be the closed normal subgroup generated by the image of a ramification map $\rho $ at $p$ (see Paragraph \ref{notrem}).  Note that since $K$ is by definition normal, the subgroup $K$ is independent of the choice of $\rho $ in $c_p$.
\end{pg}

\begin{prop}\label{P:6.2} The subgroup $K$ is the kernel of the
  surjection $\pi _1^t(C^\circ , \bar \eta )\rightarrow \pi _1^t(C, \bar
  \eta ).$
\end{prop}
\begin{proof}
  Let $U^\circ \rightarrow C^\circ $ be a finite \'etale tame cover, and
  let $U\rightarrow C$ be the normalization.  The action of $\pi
  _1^t(C^\circ , \bar \eta )$ on $U_{\bar \eta }$ factors through $\pi
  _1^t(C, \bar \eta )$ if and only if $U$ is \'etale over $p$.
This holds if and only if the pullback
$$
U\times _C\Sp (\mls O_{C, p}^{\text{sh}})\rightarrow \Sp (\mls O_{C, p}^{\text{sh}})
$$
is \'etale.  Since the formation of normalization commutes with
ind-\'etale base change, this is true if and only if the \'etale cover
$$
  U^\circ \times _{C^\circ }\Sp (\mls O_{C, p}^{\text{sh}})^\circ\rightarrow \Sp (\mls O_{C, p}^{\text{sh}})^\circ
$$
is trivial.  This in turn is equivalent to saying that $\widehat
{\Z}(1)$ acts trivially on $U_{\bar \eta }$.
\end{proof}

\begin{pg} An important property of the ramification maps is that they
  behave well in families.  Let $S= \Sp (V)$ be the spectrum of a
  discrete valuation ring with separably closed residue field, and let
  $s$ (resp.\ $\xi $) denote the closed (resp.\ generic) point of $S$.
  Let $C\rightarrow S$ be a proper smooth morphism of relative
  dimension $1$ with geometrically connected fibers.  Let $p_0, \dots
  , p_r:S\rightarrow C$ be a collection of disjoint sections, and
  write $(C_s, p_0, \dots, p_r)$ (resp.\ $(C_{\bar \xi }, p_0, \dots,
  p_r)$) for the $(r+1)$-pointed curves obtained over the geometric
  points of $S$.  Also let $C^\circ $ denote the complement $C-\{p_0,
  \dots , p_r\}$.

  Let $\bar \eta _\xi :\Sp (\Omega _\xi )\rightarrow C_\xi ^\circ $
  and $\bar \eta _s:\Sp (\Omega _s)\rightarrow C_s^\circ $ be
  geometric generic points.  By Corollary \ref{tame-inv} we then have a 
  surjection
$$
 \xymatrix{ \pi _1^t(C^\circ _{\bar \xi }, \bar \eta _\xi )\ar@{->>}[r]& \pi _1^t(C^\circ _{s}, \bar \eta _s)}
$$
well-defined up to conjugation.  For any section $p_i$, the conjugacy
class of maps $c_{p_{i, \bar \xi }}$
$$
\widehat {\Z}(1)\rightarrow \pi _1^t(C^\circ _{\bar \xi}, \bar \eta _{\xi })
$$
therefore induces a conjugacy class of maps 
\begin{equation}\label{special}
\widehat {\Z}(1)\rightarrow \pi _1^t(C_s^\circ , \bar \eta _s).
\end{equation}
\end{pg}

\begin{lem}\label{Lspecial} The conjugacy class of homomorphisms
  (\ref{special}) is equal to image of the conjugacy class of maps $c_{p_{i,
      \bar \xi }}$.
\end{lem}
\begin{proof} This is straightforward from the definition of the
  specialization maps.
\end{proof}

\section{The punctured line}\label{puncsection}

\begin{pg} Fix an isomorphism $\widehat {\Z}(1)\simeq \widehat \Z$.
  Let $\mls C$ be a class of groups satisfying the conditions in
  Paragraph \ref{2.5}.

  Let $C = \mathbb{P}^1_k$, let $\{p_1, \dots , p_r\}$ be a finite set
  of points in $\mathbb{A}^1(k)$, and let $C^\circ $ denote
  $\mathbb{A}^1-\{p_1, \dots, p_r\}$ (so $C^\circ $ is the complement
  of $r+1$ points in $\mathbb{P}^1$).  Consider the following
  properties of $\pi _1(C^\circ , \bar \eta )^{\mls C}$:
\begin{enumerate}
\item [($P_1$)] There exist embeddings $\mls O_{C,
    p_i}^{\text{sh}}\hookrightarrow \Omega $ such that the images of the
  induced ramification maps
$$
  \rho _i:\widehat {\Z}\simeq \widehat {\Z}(1)\rightarrow \pi _1(C^\circ , \bar \eta )^{\mls C}
$$
topologically generate $\pi _1(C^\circ , \bar \eta )^{\mls C}.$
\item [($P_2$)] There exist embeddings $\mls O_{C,
    p_i}^{\text{sh}}\hookrightarrow \Omega $ such that the map
$$
F_{r}\rightarrow \pi _1(C^\circ , \bar \eta )
$$
induced by the ramification maps induces an isomorphism on pro-$\mls
C$-completions.
\end{enumerate}
\end{pg}

Let $M$ be a connected locally noetherian scheme (not necessarily a
$k$-scheme), and let $p_1, \dots, p_r\in \mathbb{A}^1(M)$ be a
collection of sections such that for every (not necessarily closed)
point $m\in M$ the fibers $$p_{1, m}, \dots, p_{r, m}\in
\mathbb{A}^1(k(m))$$ are distinct.  For each $m\in M$ we can then
apply the preceding discussion to $p_{1, \bar m}, \dots, p_{r, \bar
  m}\in \mathbb{A}^1(\overline {k(m)})$.  We say that $P_1$ (resp.\
$P_2$) holds at $m$ if condition $P_1$ (resp.\ $P_2$) holds for $C_{
  \overline {k( m)}}$ and the points $\{p_{1, \bar m}, \dots, p_{r,
  \bar m}\}$.

\begin{lem}\label{cardlem} 
  (i) Let $\eta \in M$ be a point, $m\in M$ a specialization of $\eta
  $. If $P_1$ holds at $\eta $, then $P_1$ holds at $m$.

  (ii) If the orders of the elements of $\mls C$ are invertible on
  $M$, and if for some $m_0\in M$ condition $P_1$ (resp.\ $P_2$) holds
  at $m_0$, then condition $P_1$ (resp.\ $P_2$) holds at every $m\in
  M$.
\end{lem}
\begin{proof} This follows immediately from Lemma \ref{Lspecial}.
\end{proof}

\begin{prop}\label{specialprop} If $P_1$ (resp.\ $P_2$) holds for the
  set of points $\{0,1\}$ in $\mathbb{A}^1$, then $P_1$ (resp.\ $P_2$)
  holds for an arbitrary set of distinct points $\{p_1, \dots, p_r\}$ in
  $\mathbb{A}^1$.
\end{prop}
\begin{proof}
  Since the stacks $\mls M_{0, n}$ are smooth over $\Z$ it suffices by
  Lemma \ref{cardlem}(i) to consider the case when the characteristic
  of $k$ is $0$ (this assumption on the characteristic will be in
  effect until the end of this section).  Furthermore, since the
  stacks $\mls M_{0, n}$ are connected (in fact irreducible), the
  validity of $P_1$ (resp.\ $P_2$) for a set $\{p_1, \dots, p_r\}$
  depends only on the number of points $r$ and not on their specific
  values.

It follows that if $A_N$ denotes 
 $\mathbb{G}_m-\m _N(k)$, then   it
   suffices to show that if $P_1$ (resp.\ $P_2$) holds for
  $A_1$, then $P_1$ (resp.\ $P_2$) also holds for $A_N$ for any positive 
  integer $N$.

  Let $\pi _N:A_N\rightarrow A_1$ be the morphism induced by the map
  $\mathbb{P}^1\rightarrow \mathbb{P}^1$ sending $t$ to $t^N$.  The
  map $\pi _N$ is a Galois finite \'etale cover with group $\m _N(k)$.
  Let
$$
\bar \eta :\Sp (\Omega )\rightarrow A_1
$$
be a geometric generic point, and fix a compatible collection
$\{t^{1/N}\}$ of $N$-th roots of $t$ in $\Omega$.  This choice defines
a lifting
$$
\bar \eta _N:\Sp (\Omega )\rightarrow A_N
$$
of $\bar \eta $ for every $N$.  We then get a short exact sequence
$$
1\rightarrow \pi _1(A_N, \bar \eta _N)\rightarrow \pi _1(A_1, \bar
\eta )\rightarrow \m _N(k)\rightarrow 1,
$$
and hence by Lemma  \ref{L:exact} a short exact sequence
$$
  1\rightarrow \pi _1(A_N, \bar \eta _N)^{\mls C}\rightarrow \pi _1(A_1, \bar \eta )^{\mls C}\rightarrow \m _N(k)\rightarrow 1.
$$

Choose embeddings $\iota _0:\mls O_{\mathbb{P}^1,
  0}^{\text{sh}}\hookrightarrow \Omega $ and $\iota _1:\mls
O_{\mathbb{P}^1, 1}^{\text{sh}}\hookrightarrow \Omega $ defining
specialization maps
$$
\rho _0, \rho _1:\widehat \Z\rightarrow \pi _1(A_1, \bar \eta )
$$
such that the induced map $\tau :F_2\rightarrow \pi _1(A_1, \bar \eta
)$ induces a surjection (resp.\ isomorphism) on pro-$\mls
C$-completions (by assumption this is possible).  The image of $\rho _0(1)$ in $\m _N(k)$ is a generator
of $\m _N(k)$ (the generator defining the specified isomorphism
$\widehat {\Z}\simeq \widehat {\Z}(1)$).  The image $\rho _1(1)$ of
$1$ in $\m _N(k)$ is the identity.  

\begin{lem} Let $F_2$ be the free profinite group on generators
  $\gamma _0$ and $\gamma _1$, and let $\tau :F_2\rightarrow \pi
  _1(A_1, \bar \eta )$ be the map sending $\gamma _0$ to $\rho _0(1)$
  and $\gamma _1 $ to $\rho _1(1)$.  Then the kernel of the composite
  map
$$
\xymatrix{
F_2\ar[r]^-{\tau }& \pi _1(A_1, \bar \eta )\ar[r]& \m _N(k)}
$$
is 
 the topologically  free group on the generators $\gamma _0^N$ and the elements
$\gamma _0^i\gamma _1\gamma _0^{-i}$ for $0\leq i\leq N-1$.  
\end{lem}
\begin{proof}
  Let $F_2^{\text{disc}}$ denote the discrete free group on two
  generators $\gamma _0$ and $\gamma _1$ and let
  $a:F_2^{\text{disc}}\rightarrow F_2$ be the natural map.  An
  elementary computation shows that the kernel $K^{\text{disc}}$ of
  the composite
$$
F_2^{\text{disc}}\rightarrow F_2\rightarrow \m _N(k)
$$
is generated by the elements $\gamma _0^N$ and $\gamma _0^i\gamma
_1\gamma _0^{-i}$ for $0\leq i\leq N-1$.  We therefore obtain a
surjection $b:F_{N+1}^{\text{disc}}\surj K^{\text{disc}}$.  Taking the
profinite completion of the sequence
$$
F_{N+1}^{\text{disc}}\rightarrow F_2^{\text{disc}}\rightarrow \m _N(k)\rightarrow 1
$$
and using the exactness property in Proposition \ref {P:exact}
 we conclude that the map
$$
F_{N+1}\rightarrow K:=\text{Ker}(F_2\rightarrow \m _N(k))
$$
defined by the elements $\gamma _0^N$ and $\gamma _0^i\gamma _1\gamma
_0^{-i}$ for $0\leq i\leq N-1$ is a surjection.  On the other hand, by
the Nielsen-Schreier Theorem \ref{T:N-S} the group $K$ is
topologically free on $N+1$ generators, and therefore by the Hopfian
property of profinite groups \ref{Hopf} the map $F_{N+1}^{\mls
  C}\rightarrow K^{\mls C}$ is an isomorphism.
\end{proof}

We therefore obtain a commutative diagram with exact rows
$$
\xymatrix{
  1\ar[r] & F_{N+1}^{\mls C}\ar[d]^{\tau '}\ar[r] & F_2^{\mls C}\ar[d]^\tau \ar[r] & \m _N(k)\ar[r]\ar[d]^{\text{id}} & 1\\
  1\ar[r]& \pi _1(A_N, \bar \eta _N)^{\mls C}\ar[r]& \pi _1(A_1, \bar
  \eta )^{\mls C}\ar[r]& \m _N(k)\ar[r]& 1.}
$$
It follows that if $\tau $ is a surjection (resp.\ isomorphism), then
$\tau '$ is a surjection (resp.\ isomorphism).  To complete the proof
of the proposition it therefore suffices to show that the elements
$\tau (\gamma _0^N)$ and $\tau (\gamma _0^i\gamma _1\gamma _0^{-i})$
in $\pi _1(A_N, \bar \eta _N)$ are equal to the images of $1$ under
suitable choices of ramification maps.

The choice of $N$-th root $t^{1/N}$ defines a
morphism
$$
\iota _0^N:\Sp (\Omega )\rightarrow \Sp (\mls O_{\mathbb{P}^1, 0}^{\text{sh}})^\circ 
$$
such that the induced diagram
$$
  \xymatrix{& \Sp (\mls O_{\mathbb{P}^1, 0}^{\text{sh}})^\circ \ar[r]\ar[d]& A_N\ar[d]^{\pi _N}\\
    \Sp (\Omega )\ar[ru]^{\iota _0^N}\ar[r]^{\iota _0}& \Sp (\mls O_{\mathbb{P}^1, 0}^{\text{sh}})^\circ \ar[r]& A_1}
$$
commutes.  It follows that the image of $1$ under the specialization map
$$
\rho _0^N:\widehat {\Z}\rightarrow \pi _1(A_N, \bar \eta _N)
$$
induced by $\iota _0^N$ is equal to $\tau (\gamma _0)^N$.

Fix a lifting $\widetilde\gamma _0\in \text{Gal}(\Omega /k(A_1))$ of
$\tau (\gamma _0)\in \pi _1(A_1, \bar \eta )$, and let $\zeta \in \m _N(k)$ be
the image of $1$ under the composite $\widehat {\Z}\rightarrow
\widehat {\Z}(1)\rightarrow \m _N(k)$.  Since the map
$\mathbb{P}^1\rightarrow \mathbb{P}^1$ sending $t$ to $t^N$ is \'etale
over $1$, there exist unique embeddings
$$
  \iota _{\zeta ^i}^N:\mls O_{\mathbb{P}^1, \zeta ^i}^{\text{sh}}\hookrightarrow \Omega , \ \ i=0, \dots, N-1
$$
such that the diagrams
\begin{equation}\label{check1}
  \xymatrix{& \Sp (\mls O_{\mathbb{P}^1, 1}^{\text{sh}})^\circ \ar[r]\ar[d]& A_N\ar[d]^{\pi _N}\\
    \Sp (\Omega )\ar[ru]^{\iota _1^N}\ar[r]^{\iota _1}& \Sp (\mls O_{\mathbb{P}^1, 1}^{\text{sh}})^\circ \ar[r]& A_1,}
\end{equation}
and
\begin{equation}\label{check2}
  \xymatrix{\Sp (\Omega )\ar[rr]^{\widetilde\gamma_0^i}\ar[d]^{\iota _1^N}&& \Sp (\Omega )\ar[d]_{\iota _{\zeta ^i} ^N}\\
    \Sp (\mls O_{\mathbb{P}^1, 1}^{\text{sh}})^\circ \ar[d]\ar[rr]^{t\mapsto \zeta ^it}&& \Sp (\mls O_{\mathbb{P}^1, \zeta ^i}^{\text{sh}})^\circ \ar[d]\\
    A_N\ar[rr]^-{t\mapsto \zeta ^it}\ar[rd]^{\pi _N}&&A_N\ar[ld]_{\pi _N}\\
    & A_1&}
\end{equation}
commute.  If 
$$
\rho _{\zeta ^i}^N:\widehat {\Z}\rightarrow \pi _1(A_N, \bar \eta _N)
$$
denote the corresponding ramification maps, then the commutativity of
(\ref{check1}) implies that $\rho _1^N(1) = \gamma _1$ and the
commutativity of (\ref{check2}) implies that for $1\leq i\leq N-1$ we
have $\rho _{\zeta ^i}^N(1) = \gamma _0^i\gamma _1\gamma _0^{-i}$.
The completes the proof of Proposition \ref{specialprop}.
\end{proof}

\begin{rem}\label{finitegenerationremark} 
  The same argument shows that if one knows that $\pi _1^{\mls C}(A_1,
  \bar \eta )$ is finitely generated but not necessarily generated by
  the images of the ramification maps (where $A_1 =
  \mathbb{G}_m-\{0,1\}$ as in the proof), then for any integer $r$ and
  collection of distinct points $\{p_1, \dots, p_r\}$ in $\mathbb{A}^1$ the
  $\mls C$-completion of the fundamental group of $\mathbb{A}^1-\{p_1,
  \dots, p_r\}$ with respect to any base point is finitely generated.
  Indeed, it suffices to consider
  the fundamental groups $\pi _1^{\mls C}(A_N, \bar \eta _N)$, and
  these sit in a short exact sequence
  $$
  1\rightarrow \pi _1^{\mls C}(A_N, \bar \eta _N)\rightarrow \pi
  _1^{\mls C}(A_1, \bar \eta )\rightarrow \m _N(k)\rightarrow 1.
  $$
  The finite generation of $\pi _1^{\mls C}(A_N, \bar \eta _N)$
  therefore follows from the finite generation of $\pi _1^{\mls
    C}(A_1, \bar \eta )$ and the Nielsen-Schreier Theorem \ref{T:N-S}.
  Similarly, if one knows that $\pi _1^{\mls C}(A_1, \bar \eta )$ is
  topologically free then it follows that for all collections of points
  $\{p_1, \dots, p_r\}$ in $\mathbb{A}^1(k)$ the group $\pi _1^{\mls
    C}(A_N, \bar \eta _N)$ is topologically free.  The issue about
  ramification maps arises when one wants to deduce results for higher
  genus curves.
\end{rem}

\begin{rem}\label{R:5.6} The assumption that the orders of the groups
  in $\mls C$ are invertible in $k$ can be weakened as follows.  If
  one assumes $P_1$ for $\mathbb{A}^1-\{0, 1\}$ over a field of
  characteristic $0$, then it follows from the argument used in the
  proof of Proposition \ref{specialprop} that for any field $k$ and
  any set of points $p_1, \dots, p_r\in \mathbb{A}^1_k$ there exists
  embeddings $\mls O_{C, p_i}^{\text{sh}}\hookrightarrow \Omega $ such
  that the induced maps
$$
\rho _i:\widehat \Z\rightarrow \pi _1^t(C^\circ , \bar \eta )^{\mls C}
$$
topologically generate $\pi _1^t(C^\circ , \bar \eta )$.
\end{rem}

\section{Proof of Theorem \ref{theorem}$(\operatorname{i})$ and Theorem \ref{theorembis}$(\operatorname{iii})$: Wingberg}
\label{sec:an-elem-appr}

Let us first prove Theorem \ref{theorem}(i). Suppose $C$ is a proper 
smooth 
connected curve of genus $g$ over an
algebraically closed field $k$. 

\begin{prop}[Wingberg]\label{P:Wingberg}
  There is a non-canonical isomorphism $\pi_1^{\ell}(C)\cong
  S_g^{\ell}$ (where $S_g$ is defined in the paragraph preceding
  Theorem \ref{theorem}).
\end{prop}
\begin{proof}
By ``spreading out'' and using \cite[X.3.9]{SGA1}, we may assume that $k$ has positive characteristic $p\neq \ell $.  Using another spreading out argument and specialization we can further reduce to the case when $k = \overline {\mathbb{F}}_p$.  In this case the result is due to  Wingberg using the classification of Demushkin groups \cite{Wingberg} .
\end{proof}

\begin{cor}
  The group $\pi_1'(C)^{\text{\rm nil}}$ is isomorphic to
  $(S_g')^{\text{\rm nil}}$.
\end{cor}
\begin{proof}
  This follows immediately from Proposition \ref{P:Wingberg} and
 Lemma  \ref{L:pro-nilpotent-dumb}.
\end{proof}

The proof of Theorem \ref{theorembis}(iii) proceeds similarly, using 
the following substitutions: one must replace the reference to 
\cite[X.3.9]{SGA1} 
with 
Corollary \ref{tame-inv}, and the reference to \cite{Wingberg} with 
 \cite[Chapter X, 10.1.2]{NSW}, where the open case over the algebraic closure of a finite field is treated.

\section{Proof of Theorem \ref{theorembis}($\text{\rm ii}$)}

\begin{pg}
  The main result of Borne and Emsalem \cite{Borne} (building on work
  of Serre \cite{Serre}) gives an isomorphism between
  $\pi_1'(\P^1\setminus\{0,1,\infty\})^{\text{\rm sol}}$ and
  $F_2^{\prime \text{\rm sol}}$.  Given any punctured curve
  $C\subset\widebar C$, choosing a generic map $\widebar C\to\P^1$ and
  removing sufficiently many points from the base yields an open
  subcurve $U\subset C$ mapping by a finite \'etale morphism to an
  open subset $V$ of $\P^1$.  By Remark \ref{finitegenerationremark}
  the group $\pi_1'(V)^{\text{\rm sol}}$ is topologically finitely
  generated, and hence $\pi _1'(U)^{\text{sol}}$ is also topologically
  finitely generated (by the Nielsen-Schreier Theorem).  Since the
  natural map $\pi_1(U)\to\pi_1(C)$ is surjective, this completes the
  proof.  (We have omitted the base points from this argument as we
  gain no traction by keeping track of them when the given data and
  desired conclusion are so crude.)
\end{pg}

\begin{rem}
  Note that the proof of Borne and Emsalem only gives an abstract
  isomorphism $\pi_1'(\P^1\setminus\{0,1,\infty\})^{\text{\rm
      sol}}\cong F_2^{\prime \text{\rm sol}}$.  In particular, their
  methods do not show that
  $\pi_1'(\P^1\setminus\{0,1,\infty\})^{\text{\rm sol}}$ is
  topologically generated by a pair of ramification elements, and thus
  we cannot apply the techniques of \S\ref{section8} and the following
  \S\ref{sec:higher-genus-curves} to deduce results for the maximal
  pro-solvable prime-to-$p$ quotient which are stronger than Theorem
  \ref{theorembis}$(\operatorname{ii})$.
\end{rem}

\section{Higher genus curves}
\label{sec:higher-genus-curves}

Next we study higher genus curves. Let $\mls C$ be a class of finite
groups satisfying the assumptions in Paragraph \ref{2.5}, and assume
further that for every $G\in \mls C$ the order of $G$ is invertible in
$k$.  Let $C/k$ be a connected proper smooth curve of genus $g$, and let $\{p_0,
\dots, p_n\}$ be a nonempty set of points of $C(k)$.  Let $C^\circ $
denote the complement $C-\{p_0, \dots, p_n\}$.  Fix a geometric
generic point
  $$
\bar \eta :\Sp (\Omega )\rightarrow C^\circ .
$$

Consider the following properties of the $(n+1)$-pointed curve $(C, p_0,
\dots, p_n)$:
\begin{enumerate}
\item [($P_3$)] There exists a homomorphism
$$
F_{2g+n}\rightarrow \pi _1(C^{\circ}, \bar \eta )
$$
inducing a surjection on pro-$\mls C$-completions.
\item [($P_4$)] There exists a homomorphism
$$
F_{2g+n}\rightarrow \pi _1(C^{\circ}, \bar \eta )
$$
inducing an isomorphism on pro-$\mls C$-completions.
\end{enumerate}

\begin{thm}\label{devissage} Suppose $P_1$ (resp.\ $P_2$) holds for
  $(\mathbb{P}^1, \{0, 1, \infty \})$.  Then for any $m\geq 0$ and any
  $(m+1)$-pointed proper smooth genus $g$ curve $(C, \{c_0, \dots, c_m\})$ property
  $P_3$ (resp.\ $P_4$) holds for $(C, c_0, \dots, c_m)$.
\end{thm}
\begin{proof}
As in the proof of Proposition \ref{specialprop} it suffices to consider the case when $k$ has characteristic $0$.

We assume $P_1$ (resp.\ $P_2$) for $\mathbb{P}^1-\{0, 1, \infty \}$, so by Proposition \ref{specialprop} we have $P_1$ (resp. $P_2$) for any dense open in $\mathbb{P}^1-\{0,1, \infty \}$.

Since the moduli stack $\mls M_{g, m+1}$ classifying smooth proper curves of genus $g$ 
with $m+1$ marked points is irreducible (Fulton has given an algebraic
proof of this fact in the appendix to \cite{Fulton}) and the \'etale
fundamental group is constant in families of pointed curves, it
suffices to exhibit for each pair $(g, m)$ a single $m+1$-marked genus
$g$ curve for which $P_3$ (resp.\ $P_4$) holds.  We do this by a
careful analysis of hyperelliptic curves, which exist in every genus
(as smooth members of the very ample linear system $|(g+1,2)|$
on $\P^1\times\P^1$, for example).
 
 Let $f:C\to\P^1$ denote a hyperelliptic curve of
genus $g$, fibered over the projective line with $2g+2$ ramification
points $\{p_0,\dots,p_{2g+1}\}\subset C$ and $2g+2$ branch points
$\{q_0,\dots,q_{2g+1}\}\subset\P^1$, with $p_i$ lying over $q_i$ for $i=0,\ldots,2g+1$.  Given a positive integer
$n$, choose $n$ distinct closed points $r_1,\ldots,r_n\in\P^1$ disjoint
from the branch points.  If $n=0$, let $\{r_i\}$ be the empty
set.  We introduce the following notations:

\begin{enumerate}
\item $L_n=\P^1\setminus\{q_{0},r_1,\ldots,r_n\}$
\item $\widetilde L_n=L_n\setminus\{q_1,\ldots,q_{2g+1}\}$
\item $C_n=C\setminus(\{p_{0}\}\cup f^{-1}\{r_1,\ldots,r_n\})$
\item $\widetilde C_n=C_n\setminus\{p_1,\ldots,p_{2g+1}\}$
\end{enumerate}

Thus $C_n\to L_n$ is a double cover with $2g+1$ ramification points
and $\widetilde C_n\to\widetilde L_n$ is a connected $\Z/2\Z$-torsor.  The curve
$C_n$ is the complement of $2n+1$ points in a hyperelliptic curve of
genus $g$.

Let $\bar \eta :\Sp (\Omega )\rightarrow C^\circ $ be a geometric
generic point, and choose maps
$$
\iota _{q_i}:\Sp (\Omega )\rightarrow \Sp (\mls O_{\mathbb{P}^1, q_i}^{\text{sh}})^\circ , \ \ i=1, \dots, 2g+1
$$
and 
$$
\iota _{r_i}:\Sp (\Omega )\rightarrow \Sp (\mls O_{\mathbb{P}^1, r_i}^{\text{sh}})^\circ , \ \ i=1, \dots, n
$$
such that the induced map
$$
\tau :F_{2g+1+n}\rightarrow \pi _1(\widetilde L_n, \bar \eta )
$$
induces a surjection (resp.\ isomorphism) on pro-$\mls C$-completions.
This is possible by Proposition \ref{specialprop}, since we are assuming that
$P_1$ (resp.\ $P_2$) holds for $\mathbb{P}^1-\{0, 1, \infty \}$.

We then obtain a commutative diagram
\begin{equation}\label{7.1.4}
  \xymatrix{
    1\ar[r] &\text{Ker}(\chi )\ar[d]^{\tau '}\ar[r]& F_{2g+1+n}\ar[d]^{\tau }\ar[r]^{\chi }& \Z/2\Z\ar[r]\ar[d]^{\text{id}}& 1\\
    1\ar[r]& \pi _1(\widetilde C_n, \bar \eta )\ar[r]& \pi _1(\widetilde L_n, \bar \eta )\ar[r]& \Z/2\Z\ar[r]& 1,}
\end{equation}
where $\chi $ is defined to be the unique homomorphism making the diagram commute.

Let $y_i\in F_{2g+1+n}$ ($i=1, \dots, 2g+1$) denote the generator
mapping to the image of $1$ in $\pi _1(\widetilde L_n, \bar \eta )$
under the ramification map defined by $\iota _{q_i}$, and let
$y_{j+2g+1}\in F_{2g+1+n}$ ($j=1, \dots, n$) denote the generator
mapping to the image of $1$ under the ramification map defined by
$\iota _{r_j}$.

\begin{lem}\label{7.2c}
  The group $\text{\rm Ker}(\chi )$ is topologically free on the
  $4g+2n+1$ generators $g_{1i}=y_1y_i$, $i=1,\ldots,2g+1$;
  $g_{j1}=y_jy_1^{-1}$, $j=2,\ldots,2g+1$; $y_i$,
  $i=2g+2,\ldots,2g+1+n$; and $y_1y_iy_1^{-1}$,
  $i=2g+2,\ldots,2g+1+n$.
\end{lem}
\begin{proof}
  By Theorem \ref{T:N-S}, we know that $\text{Ker}(\chi )$ is topologically free of
  rank $4g+2n+1$.  Since $\text{Ker}(\chi )$ is open and profinite completion
  is an exact functor, $\text{Ker}(\chi )$ is the profinite completion of its
  intersection with the discrete free group $F_{2g+1+n}^{\text{\rm
      disc}}$ generated by the ramification elements (coming from the
  fixed choices of local uniformizers).  By Theorem \ref{Hopf}, it therefore suffices to show that
  the listed elements generate the discrete kernel, where we know an
  alternative description of the group in question: an element (word) of
  $F_{2g+1+n}^{\text{\rm disc}}$ is in $\text{Ker}(\chi )$ if and only if the
  number of factors of the type $y_i^{\pm 1}$, $i=1,\ldots,2g+1$, is
  even.

  Let $s$ be an element of $\text{Ker}(\chi )\cap F_{2g+1+n}^{\text{\rm
      disc}}$.  We proceed by induction on the length $\ell (s)$ of a
  representation of $s$ as a product of the generators $y_i^{\pm 1}$,
  the case $\ell(s)=0$ being trivial; thus, assume $\ell(s)>0$.  By
  simple manipulations, each expression of the form $y_i^{\pm
    1}y_j^{\pm 1}$ with $i,j\leq 2g+1$ (and the signs in the exponents independent) is contained in the subgroup generated by the
  $g_{1i}$ and $g_{j1}$.  If the first letter of $s$ has
  the form $y_j$ for $j\geq 2g+2$, then all of the occurences of the
  $y_i$ with $i\leq 2g+1$ take place in the remaining $\ell(s)-1$
  letters, whence we are done by induction, having expressed $s$ as a
  product of the listed generators.  Suppose that the first letter is
  the form $y_i^{\pm 1}$ with $i\leq 2g+1$.  If the second letter has
  the form $y_j^{\pm 1}$ with $j\leq 2g+1$, then we are again done by
  induction.  Thus, we may assume that the second letter has the form
  $y_j^{\pm 1}$ with $j\geq 2g+2$.  We have that $s=y_i^{\pm
    1}y_j^{\pm 1}s'=(y_i^{\pm 1}y_1^{-1})(y_1y_j^{\pm
    1}y_1^{-1})(y_1s')$ with $\ell(s')=\ell(s)-2$.  Since the first
  two factors are in the subgroup generated by the listed generators
  and $\ell(y_1s')=\ell(s)-1$, we are done.
\end{proof}

\begin{lem}\label{L:nugget} The elements $g_{11}=y_1^2$ and $g_{j1}g_{1j}=y_j^2$
  ($j=2, \dots, 2g+1$) map to the identity in $\pi _1(C_n, \bar \eta
  )$ under the composite
$$
  \text{\rm Ker}(\chi )\rightarrow \pi _1(\widetilde C_n, \bar \eta )\rightarrow \pi _1(C_n, \bar \eta ).
$$
Moreover, the closed normal subgroup generated by their images in
$\pi_1(\widetilde C_n,\bar\eta)$ is precisely the kernel of
$\pi_1(\widetilde C_n,\bar\eta)\to\pi_1(C_n,\bar\eta)$.
\end{lem}
\begin{proof}
For each $i=1, \dots, 2g+1$ choose a map
$$
\iota _{p_i}:\mls O_{C, p_i}^{\text{sh}}\hookrightarrow \Omega
$$
such that the induced diagrams
$$
  \xymatrix{
    & \Sp (\mls O_{C, p_i}^{\text{sh}})^\circ \ar[r]\ar[d]&\widetilde C_n\ar[d]\\
    \Sp (\Omega )\ar[ru]^{\iota  _{p_i}}\ar[r]^{\iota _ {q_i}}& \Sp (\mls O_{\mathbb{P}^1, q_i}^{\text{sh}})^\circ \ar[r]& \widetilde L_n}
$$
commutes.  

The induced ramification map
$$
\rho _{p_i}:\widehat  \Z\rightarrow \pi _1(\widetilde C_n, \bar \eta )
$$
then sends $1$ to $g_{11}$ in the case $i=1$ and $g_{j1}g_{1j}$ in the case $j=2, \dots, 2g+2$.  Since the composites
$$
\begin{CD}
  \widehat \Z@>\rho _{p_i}>> \pi _1(\widetilde C_n, \bar \eta )@>>>
  \pi _1(C_n, \bar \eta )
\end{CD}
$$
are zero, this implies the first statement.  The second follows from
Proposition \ref{P:6.2}.
\end{proof}

Combining Lemmas \ref{7.2c} and \ref{L:nugget} we get the following.
\begin{cor} The images of $g_{j1}$ ($j=2, \dots, 2g+1$) along with
  $y_i$ and $y_1y_iy_1^{-1}$ ($i=2g+2, \dots , 2g+n+1$) topologically generate $\pi
  _1(C_n, \bar \eta )^{\mls C}$.
\end{cor}

Let $H\subset \text{Ker}(\chi )$ denote the closed normal subgroup generated
by the elements $g_{11}$ and $g_{j1}g_{1j}$ ($j=2, \dots, 2g+1$), and
let $Q$ denote the quotient $\text{Ker}(\chi )/H$.  If $\Sigma $
denotes the kernel of $\pi _1(\widetilde C_n, \bar \eta )\rightarrow
\pi _1(C_n, \bar \eta )$, then we have a commutative diagram
$$
  \xymatrix{
    1\ar[r]& H\ar[r]\ar[d]& \text{Ker}(\chi )\ar[d]\ar[r]& Q\ar[d]\ar[r]& 1\\
    1 \ar[r] & \Sigma \ar[r]& \pi _1(\widetilde C_n, \bar \eta )\ar[r]& \pi _1(C_n, \bar \eta )\ar[r]& 1,}
$$
which induces a diagram of $\mls
C$-completions
$$
  \xymatrix{
    H^{\mls C}\ar[r]\ar[d]^a& \text{Ker}(\chi )^{\mls C}\ar[d]^b\ar[r]& Q^{\mls C}\ar[d]^c\ar[r]& 1\\
    \Sigma ^{\mls C}\ar[r]& \pi _1(\widetilde C_n, \bar \eta )^{\mls C}\ar[r]& \pi _1(C_n, \bar \eta )^{\mls C}\ar[r]& 1.}
$$
We know by Lemma  \ref{L:nugget} that the closed normal subgroup of $\pi
_1(\widetilde C_n, \bar \eta )^{\mls C}$ generated by the image of
$H^{\ms C}$ equals the image of $\Sigma^{\ms C}$. It follows that if
$b$ is surjective (resp.\ an isomorphism) then $c$ is also surjective
(resp.\ an isomorphism).    
Applying Lemma \ref{L:exact} to diagram (\ref{7.1.4}) we deduce that $b$ is surjective (resp.\ an isomorphism), so $c$ is surjective (resp.\ an isomorphism).
Since (by Lemma \ref{7.2c}) $Q$ is equal to the profinite
completion of a free group on $2g+2n$ generators, and $C_n$ is the
complement of $2n+1$ points in $C$ this proves Theorem \ref{devissage} when
$m$ is even.

To prove Theorem \ref{devissage} for odd $m$, let $C'_n$ denote the union of $C_n$ and the point in the preimage of $r_1$ corresponding to the generator $y_{2g+2}$.

\begin{lem} There exists a unique point $\tilde r_1\in C(k)$ such that if 
$$
h:\Sp (\mls O^{\text{\rm sh}}_{\mathbb{P}^1, \bar r_1})\rightarrow C
$$
denotes the unique (since $C\rightarrow \mathbb{P}^1$ is \'etale over $r_1$) lifting of the map
$$
\Sp (\mls O^{\text{\rm sh}}_{\mathbb{P}^1, \bar r_1})\rightarrow \mathbb{P}^1,
$$
then the diagram
$$
\xymatrix{
\Sp (\Omega )\ar[r]^-{\bar \eta }\ar[d]^-{\iota _{r_1}}& C\ar[d]\\
\Sp (\mls O^{\text{\rm sh}}_{\mathbb{P}^1, \bar r_1})\ar[r]\ar[ru]^-{h}& \mathbb{P}^1}
$$
commutes.  Moreover, the image in $\pi _1(C_n, \bar \eta )^{\mls C}$ of $1$ under the ramification map
$$
\widehat {\Z}\rightarrow \pi _1(C_n, \bar \eta )
$$
defined by $h$ is equal to $y_{2g+2}$.
\end{lem}
\begin{proof}
The uniqueness of $\tilde r_1$ and the uniqueness of $h$ follows from the valuative criterion for separatedness applied to the morphism $C\rightarrow \mathbb{P}^1$.  The statement that the image of $1$ under the ramification map is equal to the image of $y_{2g+2}$ follows from observing that the map
$$
\pi _1(C_n, \bar \eta )^{\mls C}\rightarrow \pi _1(L_n, \bar \eta )^{\mls C}
$$
is injective, so the ramification element is determined by its image in $\pi _1(L_n, \bar \eta )^{\mls C}$.
\end{proof}

By  Proposition \ref{P:6.2}, we have an exact sequence
$$
Z\rightarrow \pi _1(C_n, \bar \eta )^{\mls C}\rightarrow \pi _1(C_n', \bar \eta )^{\mls C}\rightarrow 1,
$$
where $Z$ is the closed normal subgroup generated by $y_{2g+2}$.  Let $Q^{\prime \mls C}$ denote the free pro-$\mls C$-group generated by $g_{j1}$ ($j=2, \dots, 2g+1$), $y_i$ ($i=2g+3, \dots, 2g+n+1$), and $y_1y_iy_1^{-1}$ ($i=2g+2, \dots, 2g+n+1$), and let $Q$ be as above.  We then have a commutative diagram
$$
\xymatrix{
\widetilde Z\ar[r]\ar[d]^{d}& Q^{\mls C}\ar[d]^-c\ar[r]& Q^{\prime \mls C}\ar[r]\ar[d]^-{c'}& 1\\
Z\ar[r]&\pi _1(C_n, \bar \eta )^{\mls C}\ar[r]&\pi _1(C_n', \bar \eta )^{\mls C}\ar[r]& 1,}
$$
where $\widetilde Z\subset Q^{\mls C}$ is the closed normal subgroup generated by $y_{2g+2}$.  Since $c$ is surjective the maps $d$ and $c'$  are surjective.  It also follows that if $c$ is an isomorphism then the map $c'$ is an isomorphism as well.
\end{proof}

\begin{rem} As in Remark \ref{R:5.6}, the assumption that the groups in $\mls C$ have order invertible in $k$ can be weakened.  If one assumes that $P_1$ holds for $\mathbb{P}^1-\{0, 1, \infty \}$ over a field of characteristic $0$, then it follows that for any $(n+1)$-pointed curve $(C, p_0,
\dots, p_n)$ over an arbitrary algebraically closed field $k$ there exists a homomorphism
$
F_{2g+n}\rightarrow \pi _1(C^{\circ}, \bar \eta )
$
inducing a surjection $F_{2g+n}\rightarrow \pi _1^t(C^\circ, \bar \eta )^{\mls C}.$
\end{rem}

\section{Deduction of Theorem \ref{theorem}$(\operatorname{ii})$ from Theorem \ref{theorembis}$(\operatorname{i})$}\label{sec:bound-relat-maxim}

\begin{pg} Let $C/k$ be a proper smooth connected curve of genus $g$, and let
  $p\in C(k)$ be a point.  Denote by $C^\circ $ the complement
  $C-\{p\}$, which is an affine curve. Let
$$
K\subset \pi _1(C^\circ , \bar \eta )
$$
be the closed normal subgroup generated by the image of a ramification map at $p$.  By Proposition \ref{P:6.2} 
the subgroup $K$ is the kernel of the
  surjection $\pi _1(C^\circ , \bar \eta )\rightarrow \pi _1(C, \bar
  \eta ).$
 In particular the map on maximal pro-$\mls C$ quotients 
$$
  \pi _1(C^\circ, \bar \eta )^{\ms C}\rightarrow \pi _1(C, \bar \eta
  )^{\ms C}
$$
is surjective, with kernel generated as a closed normal subgroup by
one element.  Moreover, the induced map on abelianizations is an
isomorphism (since it is canonically identified with the natural map
on compactly supported cohomology $H^1_c(C^\circ , \widehat
\Z(1)^{\mls C})\rightarrow H^1_c(C, \widehat \Z(1)^{\mls C})$, so the
kernel must be generated by an element of the closure of the
commutator subgroup of $\pi_1(C^{\circ},\bar\eta)^{\mls C}$.  Applying
Theorem \ref{theorembis}(i) gives the first part of the statement.  To
see that the relation is nontrivial, observe that the last statement
in Paragraph \ref{2.5} implies that there exists a prime number $\ell
$ and a commutative diagram
$$
\xymatrix{
\pi_1(C^{\circ},\bar\eta)^{\mls C}\ar[d]\ar[r]& \pi _1(C, \bar \eta )^{\mls C}\ar[d]\\
\pi_1^\ell (C^{\circ},\bar\eta)\ar[r]& \pi _1^\ell (C, \bar \eta ).}
$$
To see that the relation is nontrivial, it therefore suffices to show
that it maps to a nontrivial element of $\pi_1^\ell
(C^{\circ},\bar\eta)$, which follows from Wingberg's Theorem
\ref{P:Wingberg}.
\end{pg}

\section{Deduction of Theorem \ref{higherdim} from Theorem \ref{theorembis}$(\operatorname{ii})$}\label{sec:higherdim}

\begin{pg}\label{R:strictdef}
  By definition, if $f:X\rightarrow Y$ is a smooth morphism of schemes
  and if $D\subset X$ is a divisor, we say that $D$ has \emph{relative
    normal crossings} if \'etale locally on $X$ there exists an
  \'etale morphism
$$
X\rightarrow Y\times \Sp (\Z[x_1, \dots, x_d])
$$
such that for some $s\leq d$ the divisor $D$ is equal to the inverse
image of $Z(x_1\cdots x_s)\subset Y\times \Sp (\Z[x_1, \dots, x_d]).$
We say that $D$ has \emph{strict relative normal crossings} if the
following hold:

(i) $D$ has relative normal crossings;

(ii) There is decomposition $D = \sum_iD_i$ as a locally finite sum
of effective Cartier divisors such that
for every multi-index $I$ the scheme-theoretic intersection $D_I:=
\cap _{i\in }D_i$ is smooth over $Y$.

Note that this implies that for every geometric point $\bar
y\rightarrow Y$ the divisor $D_{\bar y}\subset X_{\bar y}$ has strict
normal crossings.
\end{pg}

\begin{lem}\label{L:higherlem}
  To prove Theorem \ref{higherdim}, it suffices to assume that $X$ is the
  complement of a strict normal crossings divisor $D$ in a smooth
  projective variety $\widebar X$.  More generally, the truth of
  Theorem \ref{higherdim} descends through arbitrary alteration (in the
  sense of de Jong \cite{dejong}) of $X$.
\end{lem}
\begin{proof}
  Let $p:X'\rightarrow X$ be the result of applying Chow's lemma and
  an appropriate alteration, so that $X'$ is the complement of a
  strict normal crossings divisor in a smooth projective variety.
  Then by \cite[IX.3.3 and IX.4.7]{SGA1} the morphism $p$ is of
  effective descent for finite \'etale covers.  It follows from
  \cite[IX.5.1]{SGA1} (and its proof) that there exists an integer $r$
  and a surjection
$$
  \pi _1(X')\langle g_1, \dots, g_r\rangle \rightarrow \pi _1(X)
$$
from the group obtained from $\pi _1(X')$ by introducing $r$ new
generators.  Thus, $\pi _1'(X)^{\text{sol}}$ is a quotient of
$$
\pi _1'(X')^{\text{sol}}\langle g_1, \dots, g_r\rangle,
$$
so if $\pi _1'(X')^{\text{sol}}$ is topologically finitely generated
then the group $\pi _1'(X)^{\text{sol}}$ is also topologically finitely
generated.
\end{proof}

After possibly replacing $X$ by an alteration, we may assume that $X$ satisfies the conditions of Lemma \ref{L:higherlem}.

\begin{lem}\label{L:step1} Let $\overline V\subset \overline X$ be a dense open subset of $\overline X$.  Then Theorem \ref{higherdim} for $V:= \overline V\cap X$ implies Theorem \ref{higherdim} for $X$.
\end{lem}
\begin{proof}
  The natural map $\pi _1(V)\rightarrow \pi _1(X)$ is surjective since
  $V$ is dense in the normal $X$, and therefore the induced map on
  prime to $p$ quotients is also surjective.  It follows that if $\pi
  _1'(V)^{\text{sol}}$ is topologically finitely generated then the
  same holds for $\pi _1'(X)^{\text{sol}}$.
\end{proof}

Let $D_1, \dots, D_r$ be the irreducible components of the divisor
$D$, and for a multi-index $I$ let $D_I$ denote $\bigcap _{i\in
  I}D_i$. By convention, if $I = \emptyset $ then $D_I:= \overline X$.
Let $L$ be a very ample line bundle on $\overline X$.  After possibly
replacing $L$ by a large power of itself, we may assume that the
natural map
\begin{equation}\label{E:surjmap}
H^0(\overline X, L)\rightarrow H^0(\overline X, L\otimes \mls O_{D_I})
\end{equation}
is surjective for every multi-index $I$.  Let $\mathbb{P}$ denote the
projective space of lines in $H^0(\overline X, L)$, and for a
subscheme $Z\subset\overline X$ let
$$
\widehat Z\subset \overline X\times \mathbb{P}
$$
denote the incidence correspondence given by the set of
(scheme-valued) points $(z, \ell )$, where $z\in Z$ and $\ell $ is a
line in $H^0(\overline X, L)$ such that $x$ is contained in the
vanishing locus of $\ell $. We will write $\widehat X$ for
$\widehat{\overline X}$.  
For a multi-index $I$, define $\widehat
D_I\subset \widehat X$ to be the closed subscheme of pairs $(x, \ell
)$ as above with $x\in D_I$.

\begin{lem}\label{L:step2b} There exists a  dense open subset $U\subset \mathbb{P}$ such that the following hold:

  (i) For every multi-index $I$ (including $I=\emptyset $) the map
  $\widehat D_{I}\times_{\mathbb{P}}U\rightarrow U$ is smooth with geometrically
  connected fibers.

(ii) $\widehat D\times_{\mathbb{P}}U$ is a divisor in the smooth $U$-scheme $\widehat
X\times_{\mathbb{P}}U$ with relative strict normal crossings over $U$.
\end{lem}

\begin{proof}[Proof of \ref{L:step2b}]
Part (i) is a consequence of Bertini's theorem and the surjectivity of the maps \ref{E:surjmap}.

Part (ii) follows immediately from part (i) using the decomposition
$\widehat D\times_{\mathbb P}U=\sum_i\widehat D_i\times_{\mathbb P}U$.
\end{proof}

Let $M\subset \mathbb{P}$ be a line with $M\cap U$ nonempty.  Since
$\widehat X\times_{\P}M\to\widebar X$ is a proper birational morphism, by Lemma
\ref{L:higherlem} we can replace $(\overline X, D)$ by $(\widehat
X\times_{\P}M,
\widehat D\times_{\P}M)$.  Applying Lemma \ref{L:step1}, we can shrink the base
of the pencil $M$ to $W:=U\cap M$ and therefore assume that there
exists a proper smooth morphism $f:\overline X\rightarrow W$ from
$\overline X$ to an open subset $W\subset \mathbb{P}^1$ such that $D$
is a divisor with relative strict normal crossings over $W$.

Let $\bar x\rightarrow X$ be a geometric point with image $\bar
y\rightarrow W$. By Theorem \ref{T:tame-seq} and Paragraph \ref{P:exact} we obtain an exact
sequence
$$
\pi _1'(X_{\bar y}, \bar x)^{\text{sol}}\rightarrow \pi _1'(X, \bar x)^{\text{sol}}\rightarrow \pi _1'(W, \bar y)^{\text{sol}}\rightarrow 1.
$$
By induction and Theorem \ref{theorembis} it follows that if $\pi
_1'(X_{\bar y}, \bar x)^{\text{sol}}$ is topologically finitely
generated then so is $\pi _1'(X, \bar x)^{\text{sol}}$.  By induction
on the dimension of $X$ we obtain Theorem \ref{higherdim}. \qed

\begin{appendix}

\section{A stacky approach to the tame fundamental group}

In this appendix we reprove some standard results on the tame
fundamental group, used in the text, by stack-theoretic methods. The main idea is to reduce various statements about open varieties to statements about proper stacks (using Abhyankar's lemma).

\begin{pg}\label{firstparab}
  For an integral separated scheme $Z$ let $\text{Fet}(Z)$ denote the category
  of finite \'etale coverings of $Z$. Let $\text{Fet}'(Z)$ denote the full
  subcategory of $\text{Fet}(Z)$ of covers whose Galois closures have degree invertible on $Z$. Both $\text{Fet}'(Z)$ and $\text{Fet}(Z)$ are Galois categories.  If $\bar z\rightarrow Z$ is a geometric point, then the inclusion $\text{Fet}'(Z)\rightarrow \text{Fet}(Z)$ corresponds under Galois duality to a surjective morphism of profinite groups $\pi _1(Z, \bar z)\rightarrow \pi _1'(Z, \bar z)$.  If $\mls C$ denotes the category of finite groups of order invertible on $Z$, then $\pi _1'(Z, \bar z)$ is the pro-$\mls C$-completion of $\pi _1(Z, \bar z)$ (in the sense of section \ref{sec:mc-c-groups}).
  
  Recall that a place of a field $\Omega $ is an equivalence class of discrete valuations on $\Omega $, where two discrete valuations $\nu _1, \nu _2:\Omega -\{0\}\rightarrow \Z$ are equivalent if they define the same topology on $\Omega $.  A valuation is \emph{trivial\/} if its image in $\Z$ is $\{0\}$, and non-trivial otherwise.  The valuation ring $R=\{\omega\in\Omega | \nu(\omega)\geq 0\}$ is a discrete valuation ring if and only if $\nu$ is non-trivial. 
  If $Z$ is an integral separated scheme with function field $\Omega $ and $\nu $ is a place of $\Omega $ defining a valuation ring $R\subset \Omega $, and if there exists a dotted arrow (necessarily unique since $Z$ is separated) filling in the diagram,
$$
\xymatrix{
\Sp (\Omega )\ar@{^{(}->}[r]\ar[rd]& \Sp (R)\ar@{-->}[d]\\
& Z,}
$$
then we call the image in $Z$ of the closed point of $\Sp (R)$ the \emph{center} of $\nu $.  In this case, a valuation on $\Omega$ is trivial if and only if its center is the generic point of $Z$.
  
  Given a set of non-trivial places  $E$ of the function field of $Z$, we let
  $\text{Fet}_E^t(Z)$ denote the full subcategory consisting of covers
  which are tamely ramified at the places $E$.  When $E\subset Z$ is a
  sum of prime divisors $E=\sum E_i$ generically lying in the normal
  locus of $Z$, we will write $\text{Fet}^t_E(Z)$ for
  $\text{Fet}^t_{\{E_i\}}(Z)$.    Given a basepoint $\ast\to Z$,
  the fundamental group of $\text{Fet}^t_E(Z)$ at $\ast$ will be
  denoted $(\pi_1)^t_E(Z,\ast)$.  
\end{pg}

\begin{pg}\label{firstpara} Let $V$ be a complete discrete valuation ring with
  separably closed residue field of characteristic $p$, $X/V$ a proper scheme with integral
  geometric fibers, and let $D\subset X$ be a divisor with relative
  normal crossings over $V$ (so that $X$ is smooth over $V$ in a neighborhood of
  $D$).  Let $\eta $ (resp $s$) denote the generic point (resp.\ closed
  point) of $S:= \Sp (V)$, and let $j:X_\eta \hookrightarrow X$ and
  $i:X_s\hookrightarrow X$ be the inclusions. Let $X^{\circ}$ denote
  $X-D$.  Let $E^{\text{abs}}$ denote the set of places of the function field of $X^\circ $ such that the induced place on the fraction field  $k(\eta )$ of $V$ has a center on $\Sp (V)$ (the ``absolute case'').  Note that any such valuation on the function field of $X^\circ $ has a center on $X$ since $X\rightarrow \Sp (V)$ is proper. In what follows we write $\text{Fet}^t(X^\circ )$ (resp. $\pi _1^t(X^\circ , \bar x)$) for $\text{Fet}^t_{E^{\text{abs}}}(X^\circ )$ (resp. $(\pi _1)_{E^{\text{abs}}}^t(X^\circ , \bar x)$).

  We will implicitly use the following throughout this appendix.
\end{pg}
  \begin{lem}\label{L:A2}
    Given $X/V$ and $D$ as above, the inclusion
    $\text{\rm Fet}^t(X^{\circ})\to\text{\rm Fet}^t_D(X^{\circ})$ is an
    equivalence of categories, and similarly for the geometric fibers
    of $X$ over $V$.
  \end{lem}
  \begin{proof}
    It suffices to prove that any covering which is tamely ramified
    along $D$ must be tamely ramified at any place $\nu $  in $E^{\text{abs}}$.   In what follows we fix a discrete valuation on $k(X )$ representing $\nu $, which we again denote (abusively) by $\nu $.  Let $x\in X$ be the center of $\nu $. If $x\notin D$ the result is trivial, so assume $x\in D$.
 Completing at a closed
    specialization of the center point (and noting that $X$ is regular
    near $D$), we are reduced to proving the analogous statement of \ref{L:A2} in  the
    case where $X$ is local and regular and $D$ has strict normal
    crossings and hence is cut out by part of a system of parameters
    for $X$, say $t_1,\ldots,t_s$.  Finally, since the residue field
    of $V$ is separably closed, we know that $V$ contains all of the
    prime-to-$p$ roots of unity. 
    
    Suppose $X=\spec A$ and let $K$ be the function field of $X$.  We
    need to show that if $K'/K$ is a finite separable extension which
    is tamely ramified along $D$, then $K'/K$ is also tamely ramified
    at $\nu $.  By Abhyankar's lemma \cite[XIII.3.6]{SGA1} there
    exists a positive integer $N$ relatively prime to $p$ such that
    the base change of $K'$ to
    $B=A[x_1,\ldots,x_s]/(x_1^N-t_1,\ldots,x_s^N-t_s)$ extends to a
    finite \'etale $B$-algebra.  Write $K''$ for the fraction field of
    $B$.  Since $K'$ extends to a finite \'etale $B$-algebra, to prove
    that $K'/K$ is tamely ramified at $\nu$ it suffices to show that
    $K^{\prime \prime }/K$ is tamely ramified at $\nu $.  Since the
    compositum of extensions which are tamely ramified at $\nu$ is
    again tamely ramified at $\nu $, this reduces the proof to the
    case when $K' = K[x]/(x^N-t)$, where $N$ is an integer prime to
    $p$ and $t$ is a parameter defining an irreducible component of
    $D$.  Let $R$ be the completion of the discrete valuation ring in
    $K$ defined by $\nu $.  We then need to show that the finite
    extension
 $$
 R\rightarrow (R[x]/(x^N-t))^\sim,
 $$
 where the right side denotes the normalization of $R[x]/(x^N-t)$, has
 ramification indices prime to $p$.

 Since $R$ contains the $N$th roots of unity, it follows by elementary
 Galois theory that the integral factors of $(R[x]/(x^N-t))^\sim $ all
 have the form $R[y]/(y^d-r)$ for some divisor $d$ of $N$ and $r$ a
 uniformizer in $R$.  For such an extension, the ramification index
 must divide the degree $d$ and therefore is prime to $p$, as desired.
\end{proof}

\begin{pg}\label{secondpara}
  The discussion of Paragraph \ref{firstpara} has an analogue when the
  base $V=\spec k$ is the spectrum of a separably closed field.  In particular, suppose
  $X/k$ is a proper geometrically integral scheme and $D\subset X$ is
  a strict normal crossing divisor.  Write $X^{\circ}=X-D$ and let
  $E^{\text{\rm abs}}$ for the set of places of the function field of
  $X$ for which every element of $k$ is a unit.  Write $\text{\rm
    Fet}^t(X^{\circ})$ for $\text{\rm Fet}_{E^{\text{\rm
        abs}}}^t(X^{\circ})$.  The proof of Lemma \ref{L:A2} readily
    yields a proof of the following.
    \begin{lem}\label{L:A2.1}
      With the preceding notation, the inclusion
    $\text{\rm Fet}^t(X^{\circ})\to\text{\rm Fet}^t_D(X^{\circ})$ is an
    equivalence of categories.
    \end{lem}
    In particular, if $Y/k$ is a geometrically integral variety with a
    geometric point $\widebar y\to Y$ and $Y\hookrightarrow X_1$ and $Y\hookrightarrow X_2$
    are two open immersions into proper geometrically integral
    $k$-varieties such that $D_i=X_i-Y$, $i=1,2$, is a strict normal
    crossings divisor, then there natural isomorphisms 
$$(\pi_1)^t_{D_1}(X_1,\widebar
y)\stackrel{\sim}{\leftarrow}\pi_1^t(Y,\widebar
y)\simto(\pi_1)^t_{D_2}(X_2,\widebar y),$$ which shows that the
computation of the tame fundamental group of $Y$ using a
compactification with simple normal crossing boundary is independent of the
compactification used.
\end{pg}

\begin{pg} The main tool in our stacky approach to the fundamental
  group is the so-called ``$N$-th root construction.''  Let $Z$ be a
  scheme and $N\geq 1$ an integer.  Let $(\mls L_1, s_1), \dots, (\mls
  L_r, s_r)$ be a collection of pairs $(\mls L_i, s_i)$, where $\mls
  L_i$ is a line bundle on $Z$ and $s_i:\mls L_i\rightarrow \mls O_Z$
  is a morphism of line bundles.  Associated to this datum is a stack
  $\mls Z_N\rightarrow Z$ defined as follows (to be more precise the
  collection $\{(\mls L_i, s_i)\}$ should also be included in the
  notation, but in all cases below there should be no ambiguity as to
  the line bundles in question).

  Consider the stack $[\mathbb{A}^1/\mathbb{G}_m]$ where
  $\mathbb{G}_m$ acts in the usual way by multiplication on
  $\mathbb{A}^1$.  It is easy to see that
  $[\mathbb{A}^1/\mathbb{G}_m]$ is isomorphic to the stack which to
  any scheme $T$ associates the groupoid of pairs $(\mls L, \alpha )$,
  where $\mls L$ is a line bundle on $T$ and $\alpha :\mls L\rightarrow
  \mls O_T$ is a morphism of $\mls O_T$--modules.  For any integer $N$
  let
$$
  p_N:[\mathbb{A}^1/\mathbb{G}_m]\rightarrow [\mathbb{A}^1/\mathbb{G}_m]
$$
be the map induced by the morphisms
$$
  \mathbb{A}^1\rightarrow \mathbb{A}^1, \ \ t\mapsto t^N, \ \ \mathbb{G}_m\rightarrow \mathbb{G}_m, \ \ u\mapsto u^N.
$$
By the definition for two integers $N$ and $M$ the map $p_{NM}$ factors as
$$
  \begin{CD}
    [\mathbb{A}^1/\mathbb{G}_m]@>p_N>> [\mathbb{A}^1/\mathbb{G}_m]@>p_M>> [\mathbb{A}^1/\mathbb{G}_m].
  \end{CD}
$$

The maps $s_i:\mls L_i\rightarrow  \mls O_Z$ define a morphism
$$
\begin{CD}
  Z@>>> [\mathbb{A}^1/\mathbb{G}_m]\times \cdots \times
  [\mathbb{A}^1/\mathbb{G}_m],
\end{CD}
$$
and we define $\mls Z_N$ to be the fiber product of the diagram
$$
\begin{CD}
  @. [\mathbb{A}^1/\mathbb{G}_m]\times \cdots \times [\mathbb{A}^1/\mathbb{G}_m]\\
  @. @VVp_N\times  p_N\times \cdots \times p_NV \\
  Z@>>>[\mathbb{A}^1/\mathbb{G}_m]\times \cdots \times
  [\mathbb{A}^1/\mathbb{G}_m].
\end{CD}
$$
If locally we trivialize the $\mls L_i$ so that the maps $s_i$ are given by sections $t_i\in \mls O_Z$, then
it follows from the definition that $\mls Z_N$ is the stack-theoretic
quotient of
$$
  \Sp (\mls O_Z[w_1, \dots, w_r]/(w_1^N=t_1, \dots , w_r^N=t_r))
$$
by the action of $\m _N\times \cdots \times \m _N$ given by
$$
(\zeta _1, \dots, \zeta_r)*w_i = \zeta _iw_i.
$$
Note the following properties of this construction:

(i) If $N$ is invertible in $Z$ then $\mls Z_N$ is a Deligne-Mumford stack.

(ii) The natural  map $\pi :\mls Z_N\rightarrow Z$ is flat and
identifies $Z$ with the coarse moduli space of $\mls Z_N$.

(iii) Over $\mls Z_N$ there is a tautological collection $\{(\mls M_i, c_i)\}$ where $c_i:\mls M_i\rightarrow \mls O_{Z_N}$ is a morphism of line bundles on $\mls Z_N$ and isomorphisms $\sigma _i:\mls M_i^N\rightarrow \pi ^*\mls L_i$ such that the diagrams
$$
\xymatrix{
\mls M_i^N\ar[rr]^-{\sigma _i}\ar[rd]^{c_i}&& \pi ^*\mls L_i\ar[ld]^-{s_i}\\
& \mls O_{\mls Z_N}&}
$$
commute.
\end{pg}

\begin{thm}\label{A2} Let the notation be as in  Paragraph \ref{firstpara}, and assume in addition that $X/V$ is smooth.  Then
  the pullback functor
  $$
    j^{\ast}:\text{\rm Fet}^t(X^{\circ})\to\text{\rm Fet}^t(X^{\circ}_{\widebar{\eta}})
  $$
  is fully faithful and the functors
  $$
    i^{\ast}:\text{\rm Fet}^t(X^{\circ})\to\text{\rm Fet}^t(X^{\circ}_s)
  $$
  and
  $$
    j^{\ast}|_{\text{\rm Fet}'(X^{\circ})}:\text{\rm Fet}'(X^{\circ})\to\text{\rm Fet}'(X^{\circ}_{\widebar{\eta}})
  $$
  are equivalences.
\end{thm}

\begin{rem} This theorem follows from the general theory of the log
  \'etale fundamental group \cite[4.7 (e)]{Illusie}.  For the
  convenience of the reader we give here a slightly different argument
  which does not (explicitly) use the theory of log geometry.
\end{rem}

\begin{proof} By blowing up along an ideal supported on $D$ we can
  assume that $D$ is a divisor with strict normal crossings (see for
  example \cite[4.2.12]{K-S}).  Let $D_1, \dots, D_r$ be the
  irreducible components of $D$, and let $\mls L_i\subset \mls O_X$ be
  the ideal sheaf (a line bundle) defining $D_i$.

  For any integer $N\geq 1$ let $\mls X_N\rightarrow X$ denote the
  result of applying the $N$-th root construction applied to $X$ with
  the inclusions of line bundles $\{\mls L_i\subset \mls
  O_X\}_{i=1}^r$.  As mentioned above, the map $\mls X_N\rightarrow X$
  is flat, the restriction $\mls X_N\times _XX^\circ \rightarrow
  X^\circ $ is an isomorphism, and if $N$ is invertible in $X$ then
  $\mls X_N$ is a Deligne-Mumford stack.

  For a Deligne-Mumford stack $\mls Y$ over $\Sp (V)$, let
  $\text{Fet}(\mls Y)$ denote the category of finite \'etale
  morphisms of algebraic stacks $\mls Y'\rightarrow \mls
  Y$, and let $\text{Fet}'(\mls Y)$ denote the category of finite
  \'etale morphisms of whose Galois closures have degree
  invertible in $V$.  (We will assume by definition that finite
  morphisms of stacks are representable, and not merely quasi-finite and
  proper.)  Note that when $\mls Y$ is a scheme these
  categories coincide with the earlier defined categories as any
  Deligne-Mumford stack admitting a representable morphism to a scheme
  is naturally isomorphic to a scheme.

  Let $\mathbb{Z}'_{\geq 1}$ denote the set of positive integers
  invertible in $V$, and view $\mathbb{Z}'_{\geq 1}$ as a category in
  which there exists a unique morphism $N\rightarrow M$ if $N$ divides
  $M$ and no morphism $N\rightarrow M$ otherwise.  The $\mls X_N$
  define a functor
\begin{equation}\label{2functor}
  \mathbb{Z}'_{\geq 1}\rightarrow (\text{stacks over $X$}).
\end{equation}

\begin{rem} The reader concerned about the $2$-categorical nature of
  the right side of (\ref{2functor}) should note that for any $N$ and
  $M$ a morphism $\mls X_N\rightarrow \mls X_M$, if it exists, is
  unique up to unique isomorphism.
\end{rem}

\begin{prop}\label{P:A8} The functors
  \begin{equation}\label{func1}
    \varinjlim _{N\in \mathbb{Z}'_{\geq 1}}\text{\rm Fet}(\mls X_N)\rightarrow \text{\rm Fet}^t(X^{\circ}), 
  \end{equation}
  \begin{equation}\label{func2}
    \varinjlim _{N\in \mathbb{Z}'_{\geq 1}}\text{\rm Fet}(\mls X_{N, \bar \eta })\rightarrow \text{\rm Fet}^t(X^{\circ}_{\bar \eta}),
  \end{equation}
  and
  \begin{equation}\label{func3}
    \varinjlim _{N\in \mathbb{Z}'_{\geq 1}}\text{\rm Fet}(\mls X_{N, s})\rightarrow \text{\rm Fet}^t(X^{\circ}_s)
  \end{equation}
  are equivalences of categories.  Moreover, (\ref{func1}),
  (\ref{func2}), and (\ref{func3}) induce equivalences between the
  subcategories $\text{\rm Fet}'$.
\end{prop}
\begin{proof} We prove that (\ref{func1}) is an equivalence of
  categories leaving the proof that (\ref{func2}) and (\ref{func3}) to the
  reader (use exactly the same argument).

  For any $N\in \mathbb{Z}'_{\geq 1}$ the open immersion
  $j:X^{\circ}\hookrightarrow \mls X_N$ is dense and $\mls X_N$ is
  normal.
  It follows that for any finite \'etale morphism $Z\rightarrow \mls
  X_N$, the stack $Z$ is equal to the normalization of $Z^{\circ}:= Z\times
  _{\mls X_N}X^{\circ}$ in $\mls X_N$.  This implies that (\ref{func1}) is
  fully faithful, and in addition shows that in order to prove that
  (\ref{func1}) is an equivalence we have to show the following statement: For any object
  $Z^{\circ}\rightarrow X^{\circ}$ of $\text{Fet}^t(X^{\circ})$ there exists $N\in
  \mathbb{Z}'_{\geq 1}$ such that the normalization of $Z^{\circ}$ in $\mls
  X_N$ is \'etale over $\mls X_N$.  This follows from the relative
  Abhyankar's lemma \cite[XIII.3.6]{SGA1}.  (Note that in spite of the fact that we allow covers with degrees
  divisible by $p$, it really does suffice to consider only covers
  branched to prime-to-$p$ orders to split the ramification.  This is
  pointed out in the proof of [\emph{loc.\ cit.\/}].)
\end{proof}
From this proposition it follows that in order to prove Theorem \ref{A2} it
suffices to prove the following stack-version of \cite[X.3.9]{SGA1}.
\end{proof}

\begin{thm} Let $\mls X\rightarrow \Sp (V)$ be a proper smooth 
  Deligne--Mumford stack.  Then the functor
  $$
    j^{\ast}:\text{\rm Fet}(\mls X)\to\text{\rm Fet}(\mls X_{\widebar{\eta}})
  $$
  is fully faithful and the functors
  $$
    i^{\ast}:\text{\rm Fet}(\mls X)\to\text{\rm Fet}(\mls X_s)
  $$
  and
  $$
    j^{\ast}|_{\text{\rm Fet}'(\mls X)}:\text{\rm Fet}'(\mls
    X)\to\text{\rm Fet}'(\mls X_{\widebar{\eta}})
  $$
  are equivalences.
\end{thm}
\begin{proof} That $i^*$ is an equivalence of categories follows from
  the invariance of the \'etale site under infinitesimal thickenings
  and the Grothendieck existence theorem for stacks \cite[1.4]{Chow}.

  That $j^*|_{\text{\rm Fet}'(\mls X)}$ is an equivalence can be seen
  by the same technique used to prove \cite[X.3.8]{SGA1}.  We recall
  the argument.  If $Z_{\eta} \rightarrow \mls X_{\eta}$ is a finite
  \'etale morphism, then it follows from the purity theorem
  \cite[X.3.1]{SGA1}, the ``prime-to-$p$'' assumption, and the
  relative form of Abhyankar's lemma [\emph{loc.\ cit.\/}] that there
  is a (totally ramified) finite extension $V'\to V$ preserving the
  residue field such that the normalization of $\mls X_V$ in
  $Z_\eta|_{V'}$ is finite \'etale.  Taking the special fiber yields a
  finite \'etale covering of $\mls X_s$, which deforms to a unique
  finite \'etale covering $Z\to \mls X_{V'}$ (its essential preimage under
  $i^{\ast}$).  The restriction $Z|_{\widebar{\eta}}$ is isomorphic to
  $Z_{\eta}|_{\widebar{\eta}}$.  Thus, any element of
  $\text{Fet}'(\mls X_{\widebar{\eta}})$ arising by pullback from an
  element of $\text{Fet}'(\mls X_{{\eta}})$ is in the
  essential image of $j^{\ast}|_{\text{\rm Fet}'(\mls X)}$.  Since
  $i^{\ast}$ is an isomorphism and remains so upon any base change
  $V'\to V$, it follows that $j^{\ast}$ is essentially surjective.
  That $j^*$ (and hence $j^*|_{\text{\rm Fet}'(\mls X)}$) is fully
  faithful follows from a similar argument applied to (spaces of) morphisms
  between finite \'etale coverings of $\mls X$.
\end{proof}

\begin{cor}\label{tame-inv} If $y_\eta \in X_\eta ^{\circ}$ is a point specializing to
  $y_s\in X_s^{\circ}$, then there is a canonical surjection
  $$
    \pi _1^t(X_{\bar \eta }^{\circ}, \bar y_\eta )\surj\pi _1^t(X_{s}^{\circ}, \bar y_s)
  $$
of tame fundamental groups and a canonical isomorphism
  $$
    \pi _1'(X_{\bar \eta }^{\circ}, \bar y_\eta )\simeq \pi _1'(X_{s}^{\circ}, \bar y_s)
  $$
  of prime-to-$p$ fundamental groups.
\end{cor}

\begin{cor}\label{C:tametostack} 
  Let the notation and assumptions be as in Theorem \ref{A2}.  The
  canonical map of topological groups
\begin{equation}\label{pimap}
  \pi _1^t(X^\circ , \bar x)\rightarrow \varprojlim _{N\in \Z'}\pi _1(\mls X_N, \bar x)
\end{equation}
is an isomorphism, where the right side is given the inverse limit of
the profinite topologies.  In particular, the inverse limit topology
on the right side coincides with the profinite topology.
\end{cor}
\begin{proof}
  For $N\in \Z'$ let $\rho _N:\pi _1^t(X^\circ , \bar x)\rightarrow
  \pi _1(\mls X_N, \bar x)$ denote the surjective projection, and let
  $K_N\subset \pi _1^t(X^\circ , \bar x)$ denote the kernel.  The fact
  that the functor
$$
\varinjlim _{N\in \Z'}\text{Fet}(\mls X_N)\rightarrow \text{Fet}^t(X^\circ )
$$
is an equivalence by Proposition \ref{P:A8} implies that $\cap _{N\in
  \Z'}K_N = \{1\}$.  It follows that the map (\ref{pimap}) is
injective.  Now the inverse limit topology makes the codomain group
Hausdorff, and the termwise surjectivity makes the map in question
have dense image (with respect to the inverse limit topology).  It
follows that the map is a closed bijection and therefore a topological
isomorphism.
\end{proof}

\begin{rem}\label{R:A.12} The same argument proving Corollary
  \ref{C:tametostack} shows that if $k$ is an algebraically closed
  field and $X/k$ is a proper smooth scheme with a divisor $D\subset
  X$ with simple normal crossings, then for any geometric point $\bar
  x\rightarrow X^\circ := X- D$ the natural map of topological groups
$$
\pi _1^t(X^\circ, \bar x)\rightarrow \varprojlim _{N\in \Z'}\pi _1(\mls X_N, \bar x)
$$
is an isomorphism, where $\mls X_N$ denotes the result of applying the
$N$-th root construction using the line bundles defining the
irreducible components of $D$.
\end{rem}

\begin{thm}\label{T:tame-seq}
  Let $f:X\to Y$ be a proper smooth morphism between integral
  Noetherian schemes with $f_{\ast}\ms O_X=\ms O_Y$.  Let $D\subset X$
  be a divisor with relative normal crossings relative to $Y$, and let
  $X^{\circ}=X\setminus D$.  Choose a geometric point $\widebar x\to
  X^{\circ}$, with image $\widebar y\to Y$.  Let $D_{X_{\bar y}}:=
  D\times _Y X_{\bar y}$ be the fiber of $D$ Then there is an exact
  sequence
    $$(\pi_1)^t_{D_{X_{\widebar y}}}(X_{\widebar y}^{\circ},\widebar
    x)\to(\pi_1)^t_{D}(X^{\circ},\widebar x)\to \pi_1(Y,\widebar y)\to
    1.$$
\end{thm}
\begin{rem}
  Just as in \S{1.6} of \cite{SGA1}, one deduces the more general
  exact sequence in low-degree homotopy groups (continuing to the
  right with $\pi_0$) when the Stein factorization of $f$ is
  non-trivial.  When $X$ and $Y$ are proper smooth $k$-varieties, the
  same sequence holds with the divisor subscripts omitted.
\end{rem}

\begin{proof}[Proof of \ref{T:tame-seq}]
  Recall that a morphism of stacks $f:\mls X\rightarrow \mls Y$ is
  \emph{separable} if $f$ is flat, and if for every field valued point
  $\Sp (K)\rightarrow \mls Y$ the fiber product $\mls X\times _{\mls
    Y}\Sp (K)$ is reduced.

\begin{lem} For every $N\in \Z'$, the map $\mls X_N\rightarrow Y$ is
  separable.
\end{lem}
\begin{proof}
  Since the map $\mls X_N\rightarrow X$ is flat the composite $\mls
  X_N\rightarrow X\rightarrow Y$ is also flat.  To verify the
  separability, it therefore suffices to show that for any
  field-valued point $y:\Sp (K)\rightarrow Y$ the fiber product $\mls
  X_N\times _Y\Sp (K)$ is reduced.  This is clear because this fiber
  is isomorphic to the result of applying the $N$-th root construction
  to the divisor with simple normal crossings $D_{y}\subset X_y:=
  X\times _Y\Sp (K)$.
\end{proof}

\begin{prop}
  Let $\mf f:\mls X\to\mls Y$ be a proper separable morphism of integral
  Noetherian Deligne-Mumford stacks with $\mf f_{\ast}\ms O_{\mls
    X}=\ms O_{\mls Y}$.  Let $\widebar x\to\mls X$ be a geometric point with
  image $\widebar y\to\mls Y$.  Then there is an exact sequence
  $$\pi_1(\mls X_{\widebar y},\widebar x)\to\pi_1(\mls X,\widebar
  x)\to\pi_1(\mls Y,\widebar y)\to 1.$$
\end{prop}
\begin{proof}
The proof is formally identical to that of Th\'eor\`eme X.1.3 and
Corollaire X.1.4 of \cite{SGA1}, to which we refer the reader for the details.
\end{proof}

We use this Proposition to prove Theorem \ref{T:tame-seq} as follows.  For every $N\in \Z'$ we obtain an exact sequence
$$
(S_N): \ \ \pi _1(\mls X_{N, \bar y}, \bar x)\rightarrow \pi _1(\mls X_N, \bar x)\rightarrow \pi _1(Y, \bar y)\rightarrow 1.
$$
These sequences form a filtering inverse system of exact sequences of profinite groups: if $N|M$ then there is a natural morphism of exact sequences $(S_M)\rightarrow (S_N)$.

The following Lemma now concludes the proof of Theorem \ref {T:tame-seq}.
\end{proof}

\begin{lem} The induced sequence
$$
\xymatrix{
\varprojlim _{N\in \Z'}\pi _1(\mls X_{N, \bar y}, \bar x)\ar[r]& \varprojlim _{N\in \Z'}\pi _1(\mls X_N, \bar x)\ar[r] & \pi _1(Y, \bar y)\ar[r]& 1\\
\pi _1^t(X_{\bar y}^\circ , \bar y)\ar[u]^{\simeq }& \pi _1^t(X^\circ , \bar y)\ar[u]^\simeq &&}
$$
is exact.
\end{lem}
\begin{proof}
Let $I\subset  \varprojlim _{N\in \Z'}\pi _1(\mls X_N, \bar x)$ be the image of $\varprojlim _{N\in \Z'}\pi _1(\mls X_{N, \bar y}, \bar x)$.  This is a closed subgroup of $ \varprojlim _{N\in \Z'}\pi _1(\mls X_N, \bar x)$ by Corollary \ref{R:A.12}, so it suffices to show that $I$ is dense in the kernel of the map to $\pi _1(Y, \bar y)$, which is immediate.
\end{proof}

\end{appendix}

\end{document}